\documentclass{amsart}
\usepackage{graphics,verbatim,color,amsthm}

 \usepackage[all]{xy}  
\theoremstyle{plain}
\newtheorem{theorem}{Theorem}
\newtheorem{problem}{Problem} 
\newtheorem{proposition}{Proposition}  
\newtheorem{exer}{Exercise}
\newtheorem{lemma}{Lemma}  
\newtheorem{cor}{Corollary} 
\newtheorem{definition}{Definition}

\theoremstyle{remark}

\theoremstyle{definition}

\usepackage{url}
\makeatletter
\def\url@leostyle{%
  \@ifundefined{selectfont}{\def\UrlFont{\sf}}{\def\UrlFont{\small\ttfamily}}}
\makeatother
\def\cal{\mathcal}

\def\Ri{ Riemannian }

\def\R{\mathbb{R}}

\def\C{\mathbb{C}}
\def\E{\mathbb{E}}

\def\V{\mathbb{V}}

\def\smallskip{\par\vspace{1mm}}
\def\medskip{\par\vspace{2mm}}
\def\bigskip{\par\vspace{3mm}}

\def\fr#1#2{\frac{#1}{#2}}

\newcount\equationcount
\def\thenumber{0}
\def\eq#1{\global\advance\equationcount by 1
   \def\thenumber{\number\equationcount}
                        {$$#1\eqno(\thenumber)$$}}

\def\E{{\mathbb E}}

\def\on{ orthonormal}

\newcommand{\dd}[2]
{
{{\partial #1}   \over {\partial #2}}
}
\begin{document}

\title[triangles]{The  Three-Body Problem and the Shape sphere.}
\author{Richard Montgomery}
 
\address{Dept. of Mathematics\\ University of California, Santa Cruz\\ Santa Cruz CA}

\email{rmont@ucsc.edu}

\date{November 15, 2011 (Preliminary Version)}



\thanks{.......}

\begin{abstract}
The   three-body problem defines a dynamics on the space of 
triangles in the plane.    The  shape sphere is  the moduli space of oriented similarity classes of planar triangles
and  lies inside   {\it shape space}, a Euclidean 3-space  parametrizing  oriented congruence classes of triangles. 
We derive and investigate  the geometry and dynamics induced on these spaces by  the three-body problem.
 We present two theorems concerning the three-body problem
whose discovery was made through the  shape space perspective.  \end{abstract}

\maketitle

 \section{Introduction} 
 
 In 1667 Newton \cite{Newton}  posed  the  three-body problem.   
% The story goes (backed up by correspondences from the time)  that 
 % a special sub-case of the problem known as ``The Lunar Problem'' gave him headaches
%and forced him to postpone   publication of the Principia for years. 
 Central questions concerning
the problem remain open today( see Problem  \ref{Oldest} below) despite  penetrating work  on the problem over  the intervening  centuries by 
some of our  most celebrated mathematicians,   
including  Euler, Lagrange, Laplace, Legendre, d'Alembert, Clairaut, Delanay,  Poincare, 
Birkhoff,  Seigel, 
Kolmogorov, Arnol'd, Moser,  and Smale.

The problem, in its crudest form, asks   to  solve the  Ordinary Differential Equation [ODE] 
of 
eq (\ref{N}) below.  This ODE governs the motion of three point masses 
 attracting each other through their mutual gravitational attractions.  
The positions of the three masses form the vertices of a triangle so I can think of the 
problem as  concerning moving triangles.   According to the relativity   principle
  laid out  by Galilieo,  the laws of physics are invariant under  isometries (see equations (\ref{translations}), (\ref{rotation}), (\ref{reflection}) below
  and the exercise that follows).  Isometries
  are the congruences of Euclid.  Galilean relativity thus implies that   congruent triangles 
  with congruent velocities  will have  congruent motions under Newton's equations.   
Now according to the SSS theorem of high school geometry, two triangles are congruent
if and only if their three side lengths are equal.   This suggests the question: 
is  there  a 2nd order ODE in the three side lengths which
describes the   three-body problem? 

The answer to the question just raised is `no!' Any attempt at    such an ODE    breaks down in a the vicinity of  collinear triangles.   
Here I derive   three alternative variables, $w_1, w_2, w_3$ to     use in place
of the side lengths. In these variables there is such an ODE.   
Unlike the vector $(a,b,c)$ of triangle edge lengths, the vector $(w_1, w_2, w_3)$ is not invariant under congruence.
Rather it is only  invariant  under the  slightly stronger  equivalence relation  of  
  ``oriented congruence''.   Oriented congruence exclude reflections.  Under reflection $(w_1,w_2, w_3) \mapsto (w_1, w_2, -w_3)$. 
  Two triangles
  are ``oriented congruent'' if there is   translations and rotation which takes one to the other.     
I   define shape space to be the space of oriented congruence classes of planar triangles. 
Shape space is   homeomorphic to $\R^3$ and   is  parameterized by the vector  $(w_1, w_2, w_3)$.

I  derive 2nd order   ODEs  (eq (\ref{ELshape})) for the $w_i$ which
are equivalent to    a special  case of the three-body problem ( the zero-angular momentum three-body problem).  
I call these  ODEs   the ``reduced ODEs''.

 Although shape space is homeomorphic to $\R^3$
  it is not isometric to $\R^3$: the shape space metric is not Euclidean. 
  % induced
 % by a quotient procedure from  the kinetic energy of physics, which we view as  a Euclidean inner product on
 % the three-body  configuration space  ( $\C^3$).   
Nevertheless  the shape space metric does enjoy  spherical
  symmetric.  So at the   heart of shape space  geometry is a sphere which  I call       the   {\it shape sphere}.  Its points 
    represent  {\it oriented similarity classes} of planar triangles. (Figure \ref{shape sphere}.)  
  %We will show that the shape sphere is indeed a metric two-sphere,
  %but its radius is not $1$, but rather $1/2$.
   The main purpose  of this article is
   to describe shape space,  the shape sphere,  and their relation to the three-body problem
   and then to   illustrate how a geometric understanding of these spaces has   yielded   new insights into  this age-old 
   problem.
   % age-old three-body problem.  
  % We   write down   the promised quotient ODE as eq (). 

\section{Three body dynamics}  Three point  
        masses   $m_1, m_2, m_3$  move  in space $\R^3$. Their
          positions as a function of time $t$ are denoted   by the position vectors  $q_1 (t), q_2 (t), q_3 (t) \in \R^3$.    
The three-body equations   derived by Newton   are 
 \begin{equation}
 \label{N}
\begin{aligned}
m_1 \ddot q_1 & =&  F_{21} + F_{31} \\
m_2 \ddot q_2 & =&  F_{12} + F_{32} \\
m_3 \ddot q_3 & =&  F_{23} + F_{13} .
\end{aligned}
\end{equation}
We sometimes refer to the equations themselves as ``the three-body problem''. 
On the left hand side of these equations    the double dots mean two time derivatives: $\ddot q = \fr{d ^2 q} {d t^2}$.
On the right hand side  \begin{equation}
\label{force}
F_{ij} = G m_i m_j \fr{ q_i - q_j}{r_{ij}^3} \qquad r_{ij} = |q_i -q_j|
\end{equation}
is the force exerted by mass $i$ on mass $j$.  The constant $G$ is Newton's gravitational
constant and is physically needed to make dimensions 
match up.   Being mathematicians, we  set $G = 1$.    The $m_i$ are   positive numbers. 
Equations (\ref{N}) are a system of second order equations in  9 variables, the 9 components of 
$q(t) = (q_1 (t), q_2 (t), q_3 (t)) \in \R^3 \oplus \R^3 \oplus \R^3$. 
%The three-body problem is the study of all the solutions of system ( \ref{N})
%and how they fit together in space.   

By design,   equations  (\ref{N}) are invariant under the Galilean group
which is the group of   transformations of space-time
$\R^3 \oplus \R$
generated by
\begin{eqnarray}
(q, t)  \mapsto & (q + c, t)   &:  \text{ translations,}  \label{translations} \\
\label{rotation}(q, t) \mapsto & (Rq , t),    & :  \text{rotations, }    \\
\label{reflection}(q, t) \mapsto & (\bar q , t),    & :  \text{reflection, }   \\
(q, t) \mapsto &(q , t+t_0)    &:  \text{time translations,} \\
(q, t)\mapsto & (q + vt, t)  & \text{: boosts.}
\end{eqnarray}
In the first equation   $c \in \R^3$ is a translation  vector.  
In the second equation $R$ is  a rotation matrix:   a three-by-three real matrix
satisfying $R R^T =Id$ and $det(R) =1$.  In the third equation $q \mapsto \bar q$
is any reflection, for example if $q = (x,y,z)$ then  
$\bar q  = (x, -y, z)$ is reflection about the $xz$ plane.  
The first three transformations generate the isometries of space. 

\begin{exer}
\label{ex1}
(A).   Verify that  the ODEs  (\ref{N}) are invariant
under  translation (\ref{translations})    as follows.
 Let $F: \R^3 \to \R^3$ be a translation: $F(q) = q+c$.  
Verify that if   ${\bf q}(t)= (q_1 (t) , q_2 (t), q_3 (t))$ satisfies (\ref{N})
then so does its translation:    $F({\bf q}(t)):= (F(q_1 (t)), F(q_2 (t)), F(q_3 (t)))$.

(B).  Formulate what it means for equations (\ref{N})
to be invariant under the other    generators the
Galilean group.  Verify these  invariances. 

(C).  [Scaling].  Consider the space-time scaling transformation:
$(q, t) \mapsto (\lambda q, \lambda ^a t)$ , $\lambda > 0$ which induces the
action on curves:  $q_i (t) \mapsto \lambda q_i (\lambda ^{-a} t)$.
Prove that equation \ref{N} is invariant under this  scaling transformation
if and only if  $a = 3/2$.  Compare with Kepler's third law.

(D)[Planar sub-problem] Let $P \subset \R^3$ be a plane through the origin
 Suppose that ${\bf q}(t)$ is a solution to   (\ref{N})
 and  that at some time, say time $t =0$,   all three bodies and their velocities
 lie in $P$:   $q_i (0), \dot q_i (0)  \in P$,  $i =1,2,3$. 
Show that    $q_i (t) \in P$ for all $t$ in the domain of the solution.  
%within this plane.  

\end{exer}

\section{Complex variables and Mass metric.} 
Exercise \ref{ex1} (D) asserts that  we can restrict the  three-body problem to a plane,
thus defining  the   ``planar 3-body problem''.   Choose $xy$ axes for this  plane
$P$ and then    identify $P$ 
  with the complex number line $\C$
by sending a point $(x,y) \in P$   to  the complex number $q = x+ iy \in \C$.
{\it The big  advantage of   complex notation is that   rotations  now
corresponds to  the operation of multiplication by a complex number a  of unit modulus. }
In other words, we may replace the matrix formula (eq (\ref{rotation}) for rotation 
by
$$q \mapsto u q,  u = exp(i \theta)$$
where $u$ is a unit complex number so that $\theta$ is real.
The number $\theta$ is the radian measure of the amount of rotation.
The set $u$ of all unit complex numbers forms the circle group, denoted $S^1$.   
 
We are now in the realm of Euclidean plane geometry.  The  locations $q_i \in \C$ of the 
three masses form    the   vertices of a Euclidean triangle.  So we describe the   triangle
as  a vector ${\bf q} = (q_1, q_2, q_3) \in \C^3$.  We call the   three-dimensional complex vector space $\C^3$
the space of 
 {\it of located triangles}, or {\it configuration space}.

Introduce  the mass inner product
   \begin{equation}
     \label{mass_metric}
  \langle {\bf v}, {\bf w} \rangle =  m_1 \bar v_1 w_1 + m_2  \bar v_2   w_2 + m_3 \bar v_3 w_3
\end{equation}
    on the space $\C^3$ of located triangles so that
     \begin{equation}
     \label{kinEn}
     K(\dot {\bf q})   =  \fr{1}{2} \langle {\bf \dot q},{\bf  \dot q}  \rangle :=  \fr{1}{2}  \Sigma m_i | \dot q_i |^2 
     , \qquad  \text{  (kinetic energy)} 
  \end{equation}
is the usual kinetic energy of a motion.  Here   ${\bf \dot q} = (\dot q_1, \dot q_2, \dot q_3)  \in \C^3$ 
is the vector representing the velocities
 of  the three masses.  Also form the potential energy
 \begin{equation}
\label{potential}
V({\bf q}) =-\{ \frac{m_1 m_2}{r_{12}} + \frac{m_2 m_3}{r_{23}} + \frac{m_1 m_3}{r_{13}} \}, \qquad \text{  (negative   potential energy ) }
    \end{equation} 
    Then 
    \begin{equation} H({\bf q}, {\bf \dot q}) = K({\bf \dot q}) +V({\bf q}) , \qquad \text{   (total energy) }
    \end{equation}
    is called the   energy of a motion ${\bf q}(t)$.
    \begin{proposition}
    \label{constenergy}
    The   energy $H$
is {\bf conserved}:    $H({\bf q}(t), {\bf \dot q} (t))$ is constant
along   solutions   ${\bf q}(t)$  to eq. (\ref{N}). 
(Different solutions typically have different constant energies.)
    \end{proposition}
    
    A complex vector space such as $\C^3$   becomes a real vector space when we only allow scalar multiplication
by real scalars.  And the real part $\langle \cdot , \cdot \rangle_{\R}$ of the  Hermitian mass inner product $\langle \cdot , \cdot \rangle$on $\E$ defines a real inner product
on $\C^3$.   Now a real inner product   induces a gradient
operator   taking  smooth functions $W: \C^3 \to \R$ to smooth real vector fields
$\nabla W :\C^3\to \C^3$ according to the rule 
\begin{equation}
\label{gradoperator}
\frac{d}{d\epsilon} W(({\bf q} + \epsilon {\bf h})|_{\epsilon = 0} = \langle \nabla W ({\bf q}), {\bf h}  \rangle_{\R}.
\end{equation}
In terms of real  linear orthogonal (not neccessarily orthonormal!) coordinates 
$\xi^j$, $j =1, \ldots , 6$ for $\C^3$ the gradient $\nabla W$ is    a  variation  of
the usual   coordinate  formula  
from vector calculus. Namely  
$\nabla W _j = \frac{1}{c_j} \dd{W}{\xi^j}$ where $c_j= \langle E_j , E_j \rangle$, 
Here the linear coordinates $\xi^j$
are related to an orthogonal basis $E_j$ for $\C^3$ as per usual:  ${\bf q} = \Sigma_{j=1} ^6 \xi^j E_j$.
We will take  the $\xi^j$ to  come in pairs $(x_j, y_j)$ as per $q_j= x_j+ i y_j$ so that the $c_j$ are then equal to the
$m_j$ in pairs and we get that  the components of our gradient:  $(\nabla V) _j = \frac{1}{m_j}( \dd{V}{x_j}, \dd{V}{y_j}) = 
 \frac{1}{m_j}( \dd{V}{x_j} + i  \dd{V}{y_j})$.

\begin{exer} 
(A) Show that Newton's equations (\ref{N}) can be rewritten
\begin{equation}
\label{N2}
\ddot {\bf q} = -\nabla V ({\bf q}),
\end{equation}
%with $U$ the negative potential function (\ref{potential}) described above.

(B)  Use  (A)   to prove constancy of energy, Proposition \ref{constenergy} above. 
\end{exer}

The energy is a function on {\it phase space} where we make
\begin{definition}
The phase space of the planar three-body problem is
$\C^3 \times \C^3$. Its points are written $({\bf q}, {\bf \dot q})$
so that the first copy of $\C^3$ represents positions and the second copy of $\C^3$
represents velocities.
\end{definition} 
We describe three other basic functions on phase space.  
    The {\it moment of inertia}
     \begin{equation}
   \label{MomInertia}I ({\bf q})= \langle {\bf q}, {\bf q} \rangle =  \Sigma m_i |q_i|^2
   \end{equation}  
   measures the overall    {\it size}
of a located triangle $q \in \C^3$.    
    The angular momentum is 
   \begin{equation}
   \label{angMom}
   J  = Im(\langle {\bf q}, {\bf  \dot q} \rangle ) =m_1 q_1 \wedge \dot q_1 + m_2 q_2 \wedge \dot q_2 + m_3 q_2 \wedge \dot q_2  , \qquad \text{   ( angular momentum ) }
    \end{equation}
  In this last formula we use the notation
  \begin{equation}
  (x+ iy) \wedge (u + iv) = det
\left(
\begin{array}{cc}
 x &  y    \\
 u &v   
\end{array}
\right)
= xv - yu
\label{wedge}
\end{equation}
which is also $Im(\bar z w)$ for  $z = x + i y,  w = u+ i v \in \C$.
   This wedge operation $z, w \mapsto z \wedge w$
   is the planar version of the cross product.  If $\times$ denotes the usual cross product of vectors in
 $\R^3$ the  $ (x,y, 0) \times (u,v,0) = (0,0, z \wedge w)$ so that $J$
 is the 3rd component of the usual angular momentum of physics. 
 If ${\bf 1} = (1,1,1) \in \C^3$ is   the vector generating translations,
 ${\bf 1} = (1,1,1) \in \C^3$ then the total linear momentum is  
  \begin{equation}
  \label{linMom}
P = \langle  {\bf \dot q}, {\bf 1}  \rangle = \Sigma m_i \dot q_i
 , \qquad \text{   (linear momentum) }
    \end{equation}
   and
  \begin{equation}
  \label{cofm}
 q_{cm} =  \langle {\bf q} , {\bf 1}\rangle /  \langle {\bf 1} , {\bf 1}\rangle= (m_1 q_1 + m_2 q_2 + m_3 q_3)/(m_1 + m_2 + m_3), \qquad \text{ (center of mass) }
    \end{equation}

\begin{exer}
\label{conservationlaws}
 Let  ${\bf q} (t)$ be a solution  to (\ref{N}) and  ${\bf \dot q} (t)$ its velocity. Show that: 

(A)  the total linear momentum $P ({\bf \dot q})$ is conserved.

  (B) the  moment of inertia  $I(t) = I(q(t))$ 
evolves according to the Lagrange-Jacobi equation:
$$\ddot I = 4H + 2 U({\bf q}).$$

(C) If  $P = 0$ and if $q_{cm} (0)  = 0$ then $q_{cm} (t) = 0$ for all time

(D) With the conditions of (C) in place the  
angular momentum $J = J({\bf q}, {\bf \dot q})$ is conserved

\end{exer}

\section{The two-body limit. Kepler's problem}

Set $m_3 = 0$ or    $q_3 = \infty$.  
Either way, we  throw out  the third equation of (\ref{N}) and the variable $q_3$.    
The first two equations of \ref{N}  remain 
with   $F_{31} = F_{32} = 0$.     These  two (vector) equations are known as the ``two-body problem''.  Set
$$\lambda = q_1 -q_2 \in \C , $$
  divide the first equation of (\ref{N}) by $m_1$ and the second equation of (\ref{N}) by $m_2$
and subtract it from the first to derive the single   equation
\begin{equation}
\label{Kepler}
\ddot \lambda = - c \frac{\lambda}{|\lambda|^3}, 
\end{equation}
with $c = m_1 + m_2$.  
This   equation (for any $c > 0$)  is often  called ``Kepler's problem''
although Kepler did not write out differential equations.  Its solutions are the famous conics of Kepler's first law, parameterized
according to Kepler's second law.  
The quantity
$$E = \frac{1}{2} | \dot \lambda |^2 - \frac{c}{|\lambda|}$$
is the associated conserved energy.
The motion of $\lambda$ is periodic and bounded if and only
if $E < 0$.  In this case the motion is a  circle or ellipse with focus at $\lambda = 0$
(or,  in the case of collisional motion a degerate   ellipse consisting of a line segment
with one endpoint at $\lambda =0$).  As long as $|q_3| >> |q_1|, |q_2|$ then the motion of 1 and 2  looks approximately like a two-body motion
over finite time intervals.

\section{Solutions of Lagrange and Jacobi}

\begin{figure}[h]
\scalebox{0.4}{\includegraphics{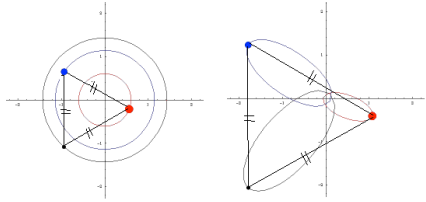}}
\caption{Lagrange solutions: the three bodies form an equilateral triangle at each instant. In
the first figure the triangle rotates about its center of mass so each individual orbit is a circle.
In the second the bodies travel on homothetic ellipses.}  \label{LagrangeOrb2}
\end{figure}

 Nearly  250 years ago 
Euler \cite{Euler} and then  Lagrange \cite{Lagrange}  wrote down explicit  solutions to the three-body problem. 
Lagrange's  solutions  are depicted in figure \ref{LagrangeOrb2}.  
These solutions of Euler and Lagrange  are the only solutions for which we have explicit analytic forms. 
The  solutions are expressed in terms of    Kepler's problem immediately above.  
  
To describe their solutions observe that we can rotate and scale a triangle ${\bf q}$
by 
\[
{\bf q} \mapsto \lambda {\bf q} , \lambda \in \C^* : = \C \setminus{0} =\text{nonzero complex numbers}.
\]
The magnitude $|\lambda|$ is the amount by which the triangle is scaled,
while the argument $\theta = Arg(\lambda)$ is the amount by which the triangle is rotated.
Now make the an\"atz  (= guess) that the  solution evolves solely by rotation and scaling:
%that a solution {\bf q}(t)  evolves by similarities:
\begin{equation}
\label{ansatz}
{\bf q}(t) = \lambda (t)  {\bf q}_0,  {\bf q}_0 \ne 0,   \lambda (t) \in \C^* := \C \setminus \{0\}
\end{equation}
\begin{exer}
\label{EulerLagrangesolutions}
Show that  ans\"atz (\ref{ansatz}) solves Newton's equations  (\ref{N}) if
and only ${\bf q}_{0, cm} : = \langle {\bf q}_0, {\bf 1} \rangle/\Sigma m_i  = 0$ , ${\bf q}_0 \ne 0$ and 
$({\bf q_0}, \lambda (t))$ solve the two equations: 
\begin{equation}
\label{Kepler}
\ddot \lambda = - c \frac{\lambda}{|\lambda|^3}, 
\end{equation}
\begin{equation}
\label{LagMult}
\frac{c}{2} \nabla I  ({\bf q}_0) = \nabla V ({\bf q}_0)
\end{equation}
where $c  = -V({\bf q}_0)/ I({\bf q}_0)$.

Hint: 
Use   form (\ref{N2}) of Newton's equation, 
the equivariance identity  $\nabla V (\lambda {\bf q}) = \frac{\lambda}{|\lambda|^3} \nabla V ({\bf q})$,  
and  Euler's identity for  homogeneous functions.
%, and the fact that $\nabla I ({\bf q} = 2 {\bf q}$. 
% $$\langle {\bf q},  \nabla U\ ({\bf q}) \rangle = - V({\bf q})$$
\end{exer} 

Eq (\ref{Kepler}) is the `Kepler problem' of the previous section. 
The second equation (\ref{LagMult}) is a Lagrange multiplier type equation  asserting
that ${\bf q}_0$ is a critical point of the function $U$, constrained to the sphere
$I = I({\bf q}_0)$.   Modulo rotations and translations, there
are exactly 5 such critical points, corresponding to the 3 collinear solutions found by Euler
and the 2 equilateral solutions of Lagrange.  They are represented by
5 points on the shape sphere.  See figure  \ref{shape sphere}.
 
%This is a good place to 
%switch gears and describe the shape space and the shape sphere.  This description
%will  allows us
%to simply describe the   5  families of solutions of Euler and Lagrange as 5 points on the shape sphere.
%See par. F of theorem 1 in the next section.
 
 %  The constant $c$ is the Lagrange multiplier for this constrained variational principle. 

 \section{Boundedness and the oldest problem}   
 
 Let us call a solution {\it bounded} if all the interparticle distances $r_{ij}$ are bounded functions of time.
 For  the two-body problem we saw  that a  solution is bounded if and only its
energy is negative.  A  partial converse of this fact   is valid in    the three-body problem. 
  
\begin{cor} (to  exercise \ref{conservationlaws}) If a  solution ${\bf q}(t)$ to eq  (\ref{N}) is bounded then its
energy $H$ is negative.
\end{cor}

Proof.  If ${\bf q}(t)$ is bounded then $I (t)$ is   bounded.
But $U > 0$ , so that if $H \ge 0$ then Lagrange-Jacobi identity $\ddot I = 4H + 2 U({\bf q})$
implies   $I(t)$ is strictly convex function and hence unbounded.  QED

\begin{figure}[h]
\scalebox{0.4}{\includegraphics{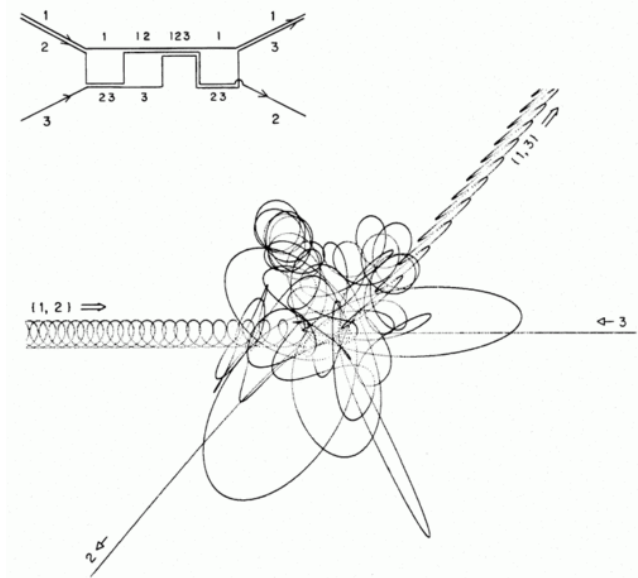}}
\caption{A typical 3 body orbit ending in a tight binary escape. Figure courtesy of Piet Hut.}  \label{Hut}
\end{figure}

The converse is false:  we know of   negative energy solutions which are unbounded.
(See figure \ref{Hut} for one.)
How false?  Many people believe ``very false''.  (See M. Hermann (\cite{Hermann} ). 

\vskip .3cm
\begin{problem}
{\bf Oldest problem in dynamical systems: }
\label{Oldest}
Are the unbounded solutions dense within the space of negative energy
solutions to the three-body problem? 
\end{problem}
\vskip .3cm

 When we say   ``the unbounded  solutions are dense''
we mean that the set of initial conditions whose solutions are unbounded forms a dense set.
This problem is completely open, despite centuries of concerted effort.
About all we know about the question  is that there is an open set of unbounded solutions,
and a set of positive measure consisting of bounded solutions.
See figure \ref{Hut} 
   for an example of what is believed to be  typical negative energy behavior
in the three-body problem: 
two   masses, say 1 and 2,  form a `tight binary' which runs away from  the third mass,
so that $r_{13}, r_{23}$   tend to infinity with time.

 \eject
 
***************************
  
PART II.  Shape Space. 

***********************************
 
%\chapter{2}

\section{Shape space. Main theorem.}

We seek a  ``reduced equations'':  a system of
three second order ODEs which encode the
  three-body  problem as a dynamical system on the space of 
congruence classes of   triangles.    The SSS theorem of elementary geometry
asserts that this     space
of congruence classes is 3-dimensional with the three  
   edge lengths  $r_{12}, r_{23}, r_{31}$ of a triangle
being coordinates.   So we expect a system of 2nd order ODEs
in the edge lengths.   However the degenerate   triangles, those with collinear vertices, 
 form a boundary for the  space of congruence classes of triangles.
 (We will see this boundary clearly in the theorem just below.)    This boundary  wreaks havoc with
dynamics.  We cannot write down   smooth   reduced differential equations
 for the dynamics of a congruence class valid  in a neighborhood of a collinear triangle.
 
 To salvage a reduced equation    I   strengthen the notion of congruence by insisting
 that the congruences be orientation preserving as well as distance preserving. Orientation preserving
 isometries are also called ``rigid motions'':  
  \begin{definition}  The group $G$ of rigid motions of the plane is the
 group of orientation preserving isometries of the plane.
 \end{definition}
\noindent Thus I exclude  reflections  from $G$. 
 Any element of $G$ is   a composition  of a rotation  and a translation.
 \begin{definition} Two planar  triangles (possibly degenerate) are `oriented congruent'
if there is a rigid motion   taking one triangle to the other. 
\end{definition} 
\begin{definition}
\noindent Shape space is the space of oriented congruence classes of triangles,
endowed with the 
%quotient topology and 
quotient metric.
\end{definition}
 \noindent With this strengthening of `congruence''  the boundary of the space of congruence classes disappears:   the collinear triangles 
are smooth points of shape space.  See Theorem \ref{main} below.
% the resulting quotient space, which we call shape  space,
%and we can    sail through them  without taking notice. 

Some words are in order regarding the meaning of   ``quotient metric''.
The space $\C^3$ of located triangles has a Euclidean metric defined by the mass metric
(eq (\ref{mass_metric})) has an associated norm $\| \cdot \|$ under which  the  Euclidean distance between
two located triangles ${\bf q}_1, {\bf q}_2$ is $\| {\bf q}_1 -    {\bf q}_2 \| $.
Our group   $G$    acts on  $\C^3$ by isometries relative to this distance and I denote the  
result of applying $g \in G$ to  ${\bf q} \in \C^3$ by $g {\bf q}$.
The  shape space metric $d$ on the quotient
by 
\begin{equation}
d([{\bf q}_1] , [{\bf q}_2] ) = \inf_{g_1, g_2 \in G}  \| g_1 {\bf q}_1 -   g_2 {\bf q}_2) \| 
\label{dist}
\end{equation}
Here    $[{\bf q}_i] $ are the `shapes', or oriented congruence classes
of the located triangles ${\bf q}_i \in \C^3$.

\begin{theorem}  (See figure \ref{shape sphere}.)
\label{main}
Shape space is homeomorphic to $\R^3$. The quotient map from the space of located triangles to
shape space  
is realized by a  map $\pi: \C^3 \to \R^3$ which is the composition
a  complex linear projection   $ \C^3 \to \C^2$ (eqs (\ref{byTranslations}))
and a real  quadratic homogeneous map $\C^2 \to \R^3$ (eq (\ref{Hopf})).  The map $\pi$  enjoys  the following
properties.  
 \begin{itemize}
\item[] A)
Two  triangles  ${\bf q}_1, {\bf q}_2 \in \C^3$ are oriented congruent  iff  $\pi({\bf q}_1) = \pi({\bf q}_2)$.
\item[]B)   $\pi$ is onto.
\item[]C)  $\pi$ projects the  triple collision locus   onto the origin.
\item[]D) If  ${\bf w} = (w_1,  w_2, w_3)$ are standard  linear coordinates on $\R^3$ then 
%the   third coordinate
$w_3$ is the signed area of the corresponding triangle, up to a mass-dependent constant.
  \item[]E)  The  collinear triangle locus  corresponds to the plane  $ w_3 = 0 $.
 \item[]F)   Let $\sigma: \R^3 \to \R^3$
be reflection across the   collinear  plane:  $ \sigma(w_1, w_2, w_3) = (w_1, w_2, -w_3)$.
Then  the two   triangles  ${\bf q}_1, {\bf q}_2 \in \C^3$ are   congruent  if and only if  either $\pi({\bf q}_1) = \pi({\bf q}_2)$
or $\pi({\bf q}_1) = \sigma(\pi({\bf q}_2))$ holds.
\item[]G)   $w_1 ^2 + w_2 ^2 + w_3 ^2 = (\frac{1}{2} I)^2$ where $I = \langle {\bf q}, {\bf q} \rangle$ (see eq (\ref{MomInertia}).
\end{itemize} 
 \end{theorem}
 
 {\it Remarks.}    E and F  of the theorem  say that the space of congruence classes of triangles
can be identified with the closed half space $w_3 \ge 0$ of $\R^3$.
The  space of collinear triangles $w_3 =0$ form its boundary, as claimed at the beginning of this section.

\subsection{The metric}

 Although shape space is homeomorphic to  $\R^3$ it is not isometric to $\R^3$.
 Shape space   geometry  is not a Euclidean geometry.    I will
 describe the shape space geometry     in more detail below in   section \ref{shapemetricsection}.
 Although the shape space geometry is not Euclidean, it is   spherically 
 symmetric, so that  the geometry of each sphere $\{R = c \}$ centered at triple collision 
 is that of the standard sphere up to a scale factor.
 I can identify that standard sphere with the {\it shape sphere}.
 
\vskip .3cm

  \subsection{The Shape Sphere}

 Add  scalings   to the group $G$  of rigid motions in order to form the 
   group of   orientation-preserving similarities whose elements are   compositions  
 of rotations, translations and scalings.   
\begin{definition}
Two planar triangles are `oriented similar' if there
is an orientation-preserving similarity  taking one to the other.
\end{definition} 
\begin{definition}
The {\it shape sphere} is the resulting quotient space of  the space of located triangles $\C^3 \setminus \C {\bf 1}$, after
the  triple collisions $\C{\bf 1}$ have been deleted.
\end{definition} 
In other words, the shape sphere is the space of oriented similarity classes of planar triangles
where I do not allow all three vertices of the triangle to coincide.

\begin{figure}[h]
\scalebox{0.4}{\includegraphics{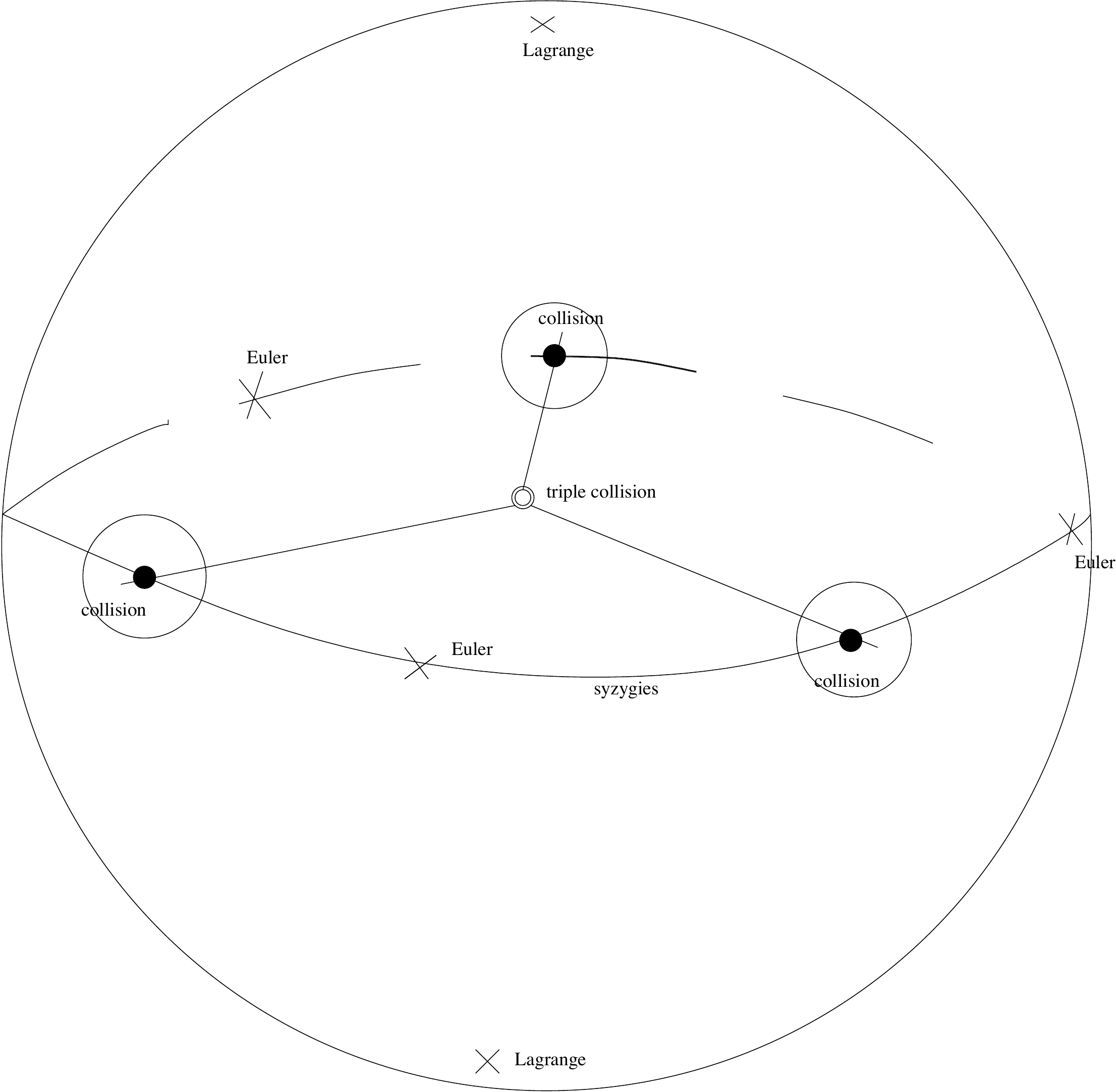}}
\caption{The shape sphere, centered on triple collision.}  \label{shape sphere}
\end{figure}

Now $\pi(\lambda {\bf q}) =   \lambda^2 \pi( {\bf q}) $ for $\lambda$ real. It follows
that   the  shape sphere can   be realized  as the space of rays through the origin
in $\R^3$. This space of rays can in turn be identified with    the  unit sphere $\|{\bf w} \| =1$   within  shape space. 
Various special types of triangles, including the 
  five families of  solutions of Euler and Lagrange are encoded on this sphere as 
  indicated in figure \ref{shape sphere}.
  
   \section{ Forming the Quotient. Proving Theorem \ref{main}.}

\begin{figure}[h]
\scalebox{0.2}{\includegraphics{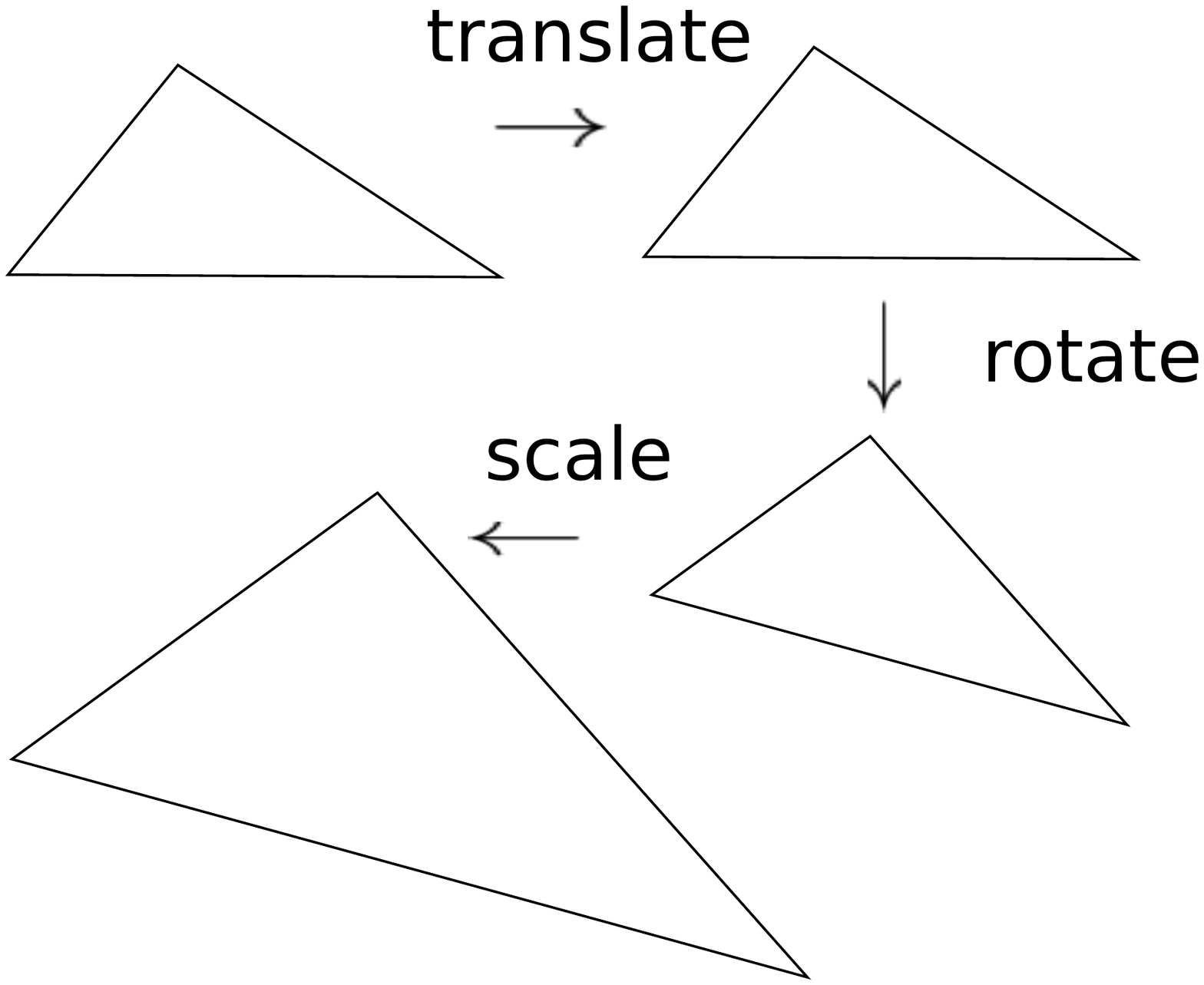}}
\caption{Translating, Rotating and Scaling a Triangle}  \label{fig_isometries}
\end{figure}

 %We begin by  showing how     the space of  (oriented) similarity classes of   triangles forms the Riemann sphere.
Recall that a vector in  $\C^3$   represent the vertices of a
planar triangle, with the ith component being the ith vertex in 
$\C \cong  \R^2$.  Translation of such a triangle
${\bf q} = (q_1, q_2, q_3) \in \C^3$  by $c \in \C$ sends  ${\bf q}$ to  the located triangle   ${\bf q} + c {\bf 1}$, 
  where ${\bf 1} = (1,1,1)$.  Rotation   by $\theta$ radians about the plane's origin
 sends   ${\bf q}$ to $e^{i \theta} {\bf q} =   ( e^{i \theta} q_1,  e^{i \theta} q_2, e^{i \theta} q_3)$.   Scaling   the plane by a positive factor $\rho$
corresponds to multiplication  by the real number $\rho$ and so sends   the triangle 
${\bf q}$ to $\rho {\bf q} =   (\rho q_1,  \rho  q_2, \rho q_3)$.  See figure \ref{fig_isometries}.

Shape space is the quotient of  $\C^3$ by the action of the group $G$ generated by
translation and rotation.  We form this quotient in two steps, translation, then rotation.  
% In the next section I  perform the  translation quotient:  
%\begin{equation}
%\pi_{tr}: \C^3 \to \C^2/ \C{\bf 1}: = \C^2 _0
%\end{equation}
% In the    section after that  we perform the quotient by rotation:  
%\begin{equation}
%\pi_{tr}: \C^2 _0 \to  \C^2 _0/ S^1 : = \R^3
%\end{equation}

\subsection{Dividing by translations}

We divide by translations by using  the isomorphism
$$\C^3/ \C {\bf 1} \cong \C{\bf 1}^{\perp}.$$
which is a special case of
$$\E/S \cong S^{\perp}$$
valid for any finite-dimensional complex vector space $\E$ with a Hermitian inner product,
and any complex linear subspace $S \subset \E$.
This isomorphism is a metric isomorphism.
Here $\E/S$ inherits a Hermitian inner product whose distance is given by 
the formula \ref{dist}
with the group  $G$  replaced by $S$  acting on $\E$ by translation, and with the elements
${\bf q}_i$ in that formula being elements of $\E$.  
In the isomorphism  the metric  I  use on $S^{\perp}$ is   the restriction of the metric from $\E$.

 In our situation $S$ is the span of ${\bf 1}$. I define   
$$\C^2 _0: =  {\bf 1} ^{\perp} = \{ {\bf q}:    m_1 q_1 + m_2 q_2 + m_3 q_3 = 0 \}$$
 the set of planar three-body configurations whose center of mass is at the origin. 
 This two-dimencional complex space   represents the quotient space
 of $\C^3$ by   translations.

 \subsection{  Jacobi coordinates: Diagonalizing the mass metric} 
 \label{Jacobi}
 
It will be   helpful
to have  coordinates diagonalizing the Hermitian form on $\C^2 _0$.
That  form is the restriction of the mass metric on $\C^3$.  Its   associated 
real positive definite quadratic form   is   the moment of inertia: 
 $I = \langle {\bf q}, {\bf q} \rangle = m_1 |q_1|^2 + m_2 |q_2 |^2 + m_3 |q_3 |^2$.
 Thus we look for coordinates $Z_1, Z_2$ on our $\C^2 _0$  such that:
   \begin{equation}
  \label{diagonalizeI}
  I = |Z_1|^2 + |Z_2|^2 \text{ whenever  } {\bf q} \perp {\bf 1}
  \end{equation}
  Jacobi found these coordinates.
 
 % (So far, we have allowed ourselves to use  any complex linear coordinates $Z_1, Z_2$  on $\C^2 _0$
%  in forming the shape coordinates $w_i$.)

  \begin{exer}
 Show that the vectors ${\bf 1} = (1,1,1),  E_1 = (\frac{1}{m_1}, -\frac{1}{m_2}, 0)$
 and $E_2 =  (\frac{-1}{m_1 +m_2}, \frac{-1}{m_1+m_2}, \frac{1}{m_3})$ form an  orthogonal 
( but not neccessarily orthonormal)  basis  relative to the mass inner product on $\C^3$. 
 \end{exer}
 The corresponding coordinates $\langle {\bf q}, {\bf 1} \rangle,  
    \langle {\bf q}, E_1 \rangle ,  \langle {\bf q} , E_2 \rangle $
    are orthogonal coordinates for $\C^3$. 
 \begin{definition}
  The  coordinates $ \langle {\bf q}, E_1 \rangle = q_1 -q_2 := Q_{12}$
 and $\langle {\bf q} , E_2 \rangle =  q_3 - \frac{m_1 q_1 + m_2 q_2}{m_1 + m_2}$
 are called   {\it  Jacobi coordinates } for $\C^2_{0}:= \{ {\bf q} \in \C^3:  {\bf q}_{cm} = 0 \}$
 relative to the partition $\{1 2 ; 3\}$ of our three masses. 
 \end{definition}
 Jacobi coordinates are indicated in  figure  \ref{Jacobi}.
  \begin{figure}[h]
\scalebox{0.5}{\includegraphics{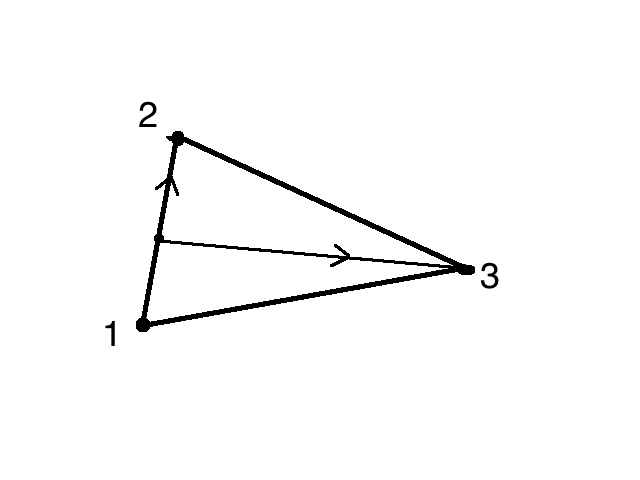}}
\caption{Jacobi vectors.}   \label{Jacobi}
\end{figure}
 
 Normalizing the  Jacobi coordinates  yields  our desired        unitary diagonalizing  coordinates
  $Z_i = \langle {\bf q} , e_i \rangle , i = 1, 2$ for $\C^2 _0$
 where   $e_i = E_i/ \| E_i\|$. 
We compute  
 \begin{equation}
 \label{Jac1}
 Z_1 = \mu_1 (q_1 -q_2)  \qquad Z_2= \mu_2 ( q_3 - \frac{m_1 q_1 + m_2 q_2}{m_1 + m_2}) 
 %= \mu_2 c_1 Q_{31} + \mu_2 c_2 Q_{32}
 \end{equation}
 with $\frac{1}{\mu_1 ^2} = \|E _1 \|^2 = \frac{1}{m_1} + \frac{1}{m_2} $
 and $\frac{1}{\mu_2 ^2} =   \|E _2 \|^2 = \frac{1}{m_3}+ \frac{1}{m_1 + m_2}$.
 %and with $c_i = m_i/ (m_1 + m_2)$.
   These normalized Jacobi coordinates  define  the complex  linear projection   
 \begin{equation}
 \label{byTranslations}
 \pi_{tr}: \C^3 \to \C^2 \qquad \pi_{tr}(q_1, q_2, q_3) = (Z_1, Z_2)
 \end{equation}
 which realizes the metric  quotient of $\C^3$ by translations.
 % and which obeys the
% quadratic relationship (\ref{diagonalizeI}).

  \subsection{Dividing by rotations}
  
  It remains to divide $\C^2 _0$ by the action of  rotations.   A rotation by $\theta$ radians acts on 
  $\C^2 _0$  by complex scalar multiplication by the unit modulus complex
  number $e^{i \theta}$.  For example, such a rotation acts on the triangle's vertices 
  $q_j$ by $q_j \mapsto e^{i \theta} q_j$ so it acts on the   triangle edges coordinates
   $Q_{jk }:= q_j - q_k $ by $Q_{jk} \mapsto e^{i \theta} Q_{jk }$.   Consequently
   it acts on the normalized Jacobi coordinates by
    $(Z_1, Z_2) \mapsto (e^{i \theta} Z_1, e^{i \theta} Z_2)$.
    
Some generality  perhaps  clarifies the situation.   Let $\V$ be a two-dimensional
 complex Hermitian space,  like our $\C^2 _0$ and  $(Z_1, Z_2)$ be Hermitian orthonormal coordinates on
 $\V$.  Then the rotation group   $S^1 =  \{e ^{ i \theta} , 0 \le \theta \le 2 \pi \}$ of  unit complex numbers
   acts on   $\V$ by scalar multiplication,  by $(Z_1, Z_2) \mapsto (e^{i \theta} Z_1, e^{i \theta} Z_2)$ as above.
 We will show that 
 the   quotient  space  $\V /S^1$   is   homeomorphic to  $\R^3$ and we will work out 
 the metric on it.    This $\R^3$ is our shape space. 
 
 Observe that the functions $Z_i \bar Z_j$, $i, j =1, 2$  remain unchanged under rotation.
 We put them together into a 2 by 2 Hermitian matrix:
 \begin{equation}
 \label{phi}
 \Phi (Z_1, Z_2) =  
\left(
\begin{array}{cc}
|Z_1|^2 & Z_1 \bar Z_2    \\
\bar Z_1 Z_2 & | Z_2|^2  
\end{array}
\right)
=   A .
\end{equation}
Or
\begin{equation}
\label{factorization}\Phi ({\bf Z})  = {\bf Z} ^t {\bf Z}^*
\end{equation}
where
 $$ {\bf Z} = (Z_1, Z_2) \qquad; {\bf Z}^t  = \left(
\begin{array}{c}
Z_1    \\
 Z_2  
\end{array} 
\right); \qquad {\bf \bar Z}  = (\bar Z_1, \bar Z_2)
$$
From the factorization (\ref{factorization}) we see:  
 $\Phi( {\bf Z}) {\bf  Z}^t = (|Z_1| ^2 + |Z_2|^2){\bf  Z}^t $
while $\Phi({\bf Z}) {\bf  W}^t  =0$ for  ${\bf W} \perp {\bf Z}$.  
Thus $\Phi(Z_1, Z_2)$ is the matrix of orthogonal projection onto  the complex  line spanned
by  ${\bf Z}$ (assuming ${\bf Z} \ne {\bf 0}$)
multiplied by $\| {\bf Z } \|^2$.   Now two nonzero vectors ${\bf Z}, {\bf U }$ are related by rotation if and only if they span the same  complex line
  and their lengths are equal.  It follows that the image of $\Phi$ is an accurate rendition of
the quotient space $\V/S^1$, with $\Phi$ being the quotient space. 
  What is the image of $\Phi$?  Well, we have just seen that it consists of the  Hermitian
  matrices of rank $1$ whose nonzero eigenvalue is positive
  (corresponding to  $\| {\bf Z } \|^2$), together with the zero matrix (corresponding to ${\bf Z} =0$).
  In terms of the determinant and trace  these conditions on $A$
  are $det(A) = 0$ and $tr(A) \ge 0$.   Let us coordinatize Hermitian matrices by  
   \begin{equation}
 \label{M} 
 A  =    
\left(
\begin{array}{cc}
w_4 + w_1  & w_2 + i w_3    \\
w_2 - i w_3 & w_4 - w_1  
\end{array}  
\right) \; , w_j \text{ real}.
\end{equation}
So that $det(A) = w_4^2 - w_1 ^2 - w_2 ^2 -w_3 ^2$
and $tr(A) = w_4$. 
  The discussion we have just had    proves:
  \begin{proposition}
  The image of the  map $\Phi$
  is the cone of  two by two Hermitian matrices $A$ 
  as above (eq \ref{M}) satisfying
  \begin{equation}
  w_4^2 - w_1 ^2 - w_2 ^2 -w_3 ^2 = 0
  \label{lightcone1}
  \end{equation}
  and 
  \begin{equation}
  w_4 \ge 0
  \label{cone2}
  \end{equation}
  This cone realizes the quotient $\V/S^1 = \C^2/S^1$
  with 
   $\Phi$ implementing the quotient map $\V \to \V/S^1$.   
  \end{proposition}
  
   Now   map the real 4 dimensional space of Hermitian matrices to $\R^3$
  onto $\R^3$ by   projecting out the trace  part  $w_4$ : 
  $$(w_1, w_2, w_3, w_4) \mapsto  pr(w_1, w_2, w_3, w_4) = (w_1, w_2, w_3).$$
  The restriction of this     projection to our   cone (eqs (\ref{lightcone1}), (\ref{cone2}) 
  is a homeomorphism onto   $\R^3$.  Indeed    solve  the cone equations for  $w_4$
  to find $w_4 = +\sqrt{ w_1 ^2 + w_2 ^2 + w_3 ^2}$ and hence
  the  
    the inverse of the restricted  projection is  
    $(w_1, w_2, w_3) \mapsto (w_1, w_2, w_3, \sqrt{ w_1 ^2 + w_2 ^2 + w_3 ^2})$
   We have proved:
    \begin{proposition} 
  \label{Hopfa} The map 
  \begin{equation}
  \label{HopfMap}
  \pi^{rot}= pr \circ \Phi: \C^2 \to \R^3
  \end{equation}
  given by 
  \begin{equation} 
 \pi^{rot} (Z_1, Z_2) =   (\frac{1}{2}(|Z_1|^2 - |Z_2|^2,  Re(Z_1 \bar Z_2), Im(Z_1 \bar Z_2) = (w_1, w_2, w_3)
  \label{Hopf} 
  \end{equation}
  realizes $\R^3$ as  the quotient space of $\C^2$ by the rotation group $S^1$.
    \end{proposition}
    
    {\bf Remark.} The restriction of the map (\ref{HopfMap}) to the sphere $w_4 =1$ is the famous {\it Hopf map}
    from the three-sphere to the two-sphere.
    
   \subsection{Proof of theorem \ref{main}.}
    We compose the projections $\pi_{tr}$ of  eq (\ref{byTranslations}) and   the map $\pi^{rot}$ immediately
    above.  The first realizes the quotient by translations and the second realizes the quotient by rotations
    so together they realize the full quotient by the group of rigid motions.    This establishes A and B of the
    theorem.   Property C, that the 
    only triangles sent to $0 \in \R^3$
    are the triple collision triangles ${\bf q} = (q, q,q)$  follows directly  from the formulae  for $\pi^{rot}$ and
    $\pi_{tr}$. Indeed, the   only  
    point of $\C^2$ mapped to $0$ by $\pi_{tr}$  is the origin $0$, and the only points of $\C^3$ mapped to the origin by
    $\pi_{tr}$ are the triple collision points. 
     
    We verify property D which says  $w_3$ is a mass-dependent constant times the
    oriented area of the triangle.  We have $w_3 = -Z_1 \wedge Z_2$.  Recall that the wedge (eq (\ref{wedge}))
    $z \wedge w = Im(\bar z w)$ represents  the oriented area of the parallelogram whose edges are$z = x+ iy$ and $w = u+ iv$.
    Thus the oriented area of our   triangle is  
    $\frac{1}{2} (Q_{21}) \wedge (Q_{31})$ where we write $Q_{ij} = q_i - q_j$ for the edge connecting  vertex j to vertex i.
    We have  $Z_1 = \mu_1 Q_{12}$ and $Z_2 = \mu_2 (p_1 Q_{31} + p_2 Q_{32})$ where $p_1 = m_1/(m_1 +m_2)$
    and $p_2 = m_2/(m_1 + m_2)$ so that $p_1 + p_2 = 1$.  Use $Q_{12} + Q_{23} + Q_{31} = 0$
    and $Q_{ij} = -Q_{ji}$ to compute that $Z_2 = \mu_2 (Q_{31} - p_2 Q_{12})$. Now the wedge operation is
    skew symmetric:  $Q_{12} \wedge Q_{12} = 0$. It follows that 
    $w_3 = - \mu_1 \mu_2 \frac{1}{2} Q_{12} \wedge Q_{31} = + \mu_1 \mu_2 \frac{1}{2} Q_{21} \wedge Q_{31}$
    as desired.    
    
    Property E follows immediately from property D.
    To establish  property F regarding the operation of reflection on triangles, observe that we can   reflect triangle ${\bf q}$
    by changing all vertices $q_i$ to $\bar q_i$ which in turn changes $(Z_1, Z_2)$
    to its conjugate vector $(\bar Z_1, \bar Z_2)$.  This conjugation operation leaves $w_1$ and $w_2$
    unchanged and changes $w_3$ to $-w_3$: the oriented area flips sign. 
    
    Property G is a computation.  Observe from eqs (\ref{phi}, \ref{M}) that $w_4 = \frac{1}{2}I$ and recall the cone condition
    eq (\ref{lightcone1}): 
    $w_4 ^2 = w_1 ^2 + w_2 ^2 + w_3 ^2$.

QED

  \section{Mechanics via Lagrangians.}
  
One of my goals here is to write down the reduced equations   encoding Newton's equations (\ref{N}) on shape space.
 My strategy for achieving this goal is   to  push the least action principle for the three-body problem down
 from the space $\C^3$ of located triangles to our  shape space $\R^3$.    
We begin   by stating the least action principle.

Any classical mechanical system can be succinctly encoded  by its  {\it Lagrangian} $L$.
(\cite{Arnold} or \cite{Feynman}), 
 \begin{equation}
 L = K -V
 \label{Lag}
 \end{equation}
 the difference of its kinetic $(K)$ and potential $(V)$ energies.
(Recall that the energy is the sum $K+V$.)  
% = \text{ kinetic } + \text{ potential}   \text{ ;  for us } V = -U)$$
Integrating the  Lagrangian over a path $c$ in the configuration space
of the mechanical system defines  that path's  {\it  action} : 
$$A[c] =  \int_c L dt  = \int_a ^b L(c (t), \dot c (t)) dt.$$
In this last expression we have taken $c$ to be parameterized by
the time interval $[a,b]$ so that  $c:[a, b] \to Q$ where  $Q$ denotes
the configuration space.  
The {\it principle of least action} asserts that a curve
satisfies Newton's equations 
if and only if $c$ minimizes  $A$
among all paths $\gamma:[a,b] \to Q$ for which
$c(a) = \gamma(a)$ and $c(b) = \gamma (b)$.

The principle is not a theorem but rather it is a guiding principle. To turn the principle  into a theorem
requires   careful wording and   more hypothesis.    Here is  such a theorem   in the
case of the three-body problem.  
\begin{theorem} 
\label{Min} If a curve $c:[0, T] \to \C^3$ minimizes the
action among all curves $\gamma: [0,T] \to \C^3$ sharing its endpoints 
and if $c$ has no collisions on the open interval $(0,T)$
then $c$ solves Newton's equations on $(0, T)$. Conversely, if  a curve $c:[a,b] \to \C^3$
satisfies Newton's equations then there is an $\epsilon > 0$
such that the restriction $c |_{[r,s]}$ of $c$ to any subinterval $[r,s]$ of size $s-r \le \epsilon$
mininimizes the action among all curves $\gamma: [r,s] \to \C^3$ sharing its endpoints:
$\gamma(r) = c(r), \gamma(s) = c(s)$.  
\end{theorem}

An analogous theorem  holds  regarding the principle of least action in situations much more general than that of the three-body problem. 
For example, we can take the configuration space    $Q = \R^n$  a real  vector space and   $K$    the squared norm
associated to any inner product  on $\R^n$.  We view $K$ as being applied to   velocities   $v \in \R^n$. Take  
       $V: \R^n \to \R$   any  smooth function.   Then the Lagrangian is   $L(x, v) = K(v) - V(x)$ which is a function
       on the   phase space  $\R^n \times \R^n$.    Newton's equations 
 are the 2nd order differential equation
  \begin{equation}
\label{Nprime}
\ddot c = - \nabla V( c).
\end{equation}
 for curves  $c$ on   $\R^n$. In these equations  
   $\nabla$ is the gradient associated to the kinetic energy $K$, 
as described in analogy to  eq (\ref{gradoperator}).   
 More generally  we can take $Q$ to be any Riemannian manifold
 with   $K$   the associated  kinetic energy ($\frac{1}{2}\Sigma g^{ij} (q)p_i p_j$.  Here the coordinate expression of 
 the   \Ri metric is $g_{ij}$.  And take  $V: Q \to \R$ any smooth function.
 Then the $\nabla$ of Newton's equations above is the covariant derivative and the second derivative of $c$
 becomes  the second covariant derivative of the curve.  
 % and the least action principle applies to Newton's equations on $Q$.
% We do not need this generality, but we will need the generality of such a $K$ as in  (4) of the theorem -- a \Ri metric on $\R^3$.

 \subsection{Euler-Lagrange Equations}

   Let $\xi^a, a = 1, \ldots , n$ be coordinates on our configuration space  $Q$.
   Then the Lagrangian is a function of the $\xi^a$ and its
   formal time derivatives $\dot \xi ^a$: 
 $$L = L (\xi^1, \ldots , \xi^n, \dot \xi ^1, \ldots , \dot \xi ^n).$$
 The  
 {\it Euler-Lagrange equations} : 
\begin{equation}
\frac{d}{dt}( \dd{L}{\dot \xi^a} )= \dd{L}{\xi^a}
\label{EL}
\end{equation}
are ODEs which a path $\xi^a = \xi^a (t)$ must   satisfy
if it minimizes the action.
They are  Newton's equations expressed in the new 
  coordinates $\xi_a$.   

  Words are in order regarding the left-hand side of the EL equations (\ref{EL}).
   We   compute  $\dd{L}{\dot \xi^a}$  by treating   $\xi^a$ and $\dot \xi^a$ as independent variables.
The resulting  $\dd{L}{\dot \xi^a}$  is now a  function of the variables  $\xi^a, \dot \xi^a$.
We then  compute $\frac{d}{dt}( \dd{L}{\dot \xi^a} )$
by   formally replacing  the independent 
variables $\xi^a, \dot \xi^a$ in $ \dd{L}{\dot \xi^a} $  by an alleged    curve $\xi^a (t)$ and its time derivatives $\dot \xi^a (t)$
so as to get a function of time which we finally   differentiate formally  using the chain rule. 
     
\begin{exer}  Suppose that $K = \frac{1}{2} \Sigma g_{a b} \dot \xi ^a \dot \xi ^b$
and that $V = V (\xi^1, \ldots , \xi^n)$.Verify that Newton's equations (\ref{N2})
are equivalent to the {\it Euler-Lagrange equations} with respect to the coordinates
$\xi^a$.
%:  $\Sigma_b g_{ab} \ddot \xi ^b =\dd{V}{\xi^a}$
\end{exer}

One of the beauties and the powers of the action principle is it is coordinate-independent.
 If a path  minimizes the  action then it does not  matter what 
 coordinate system we use to express that path.  The path  still minimizes the action
 and so satisfies the Euler-Lagrange equations in that coordinate system. 
% Regarding coordinate independence, we
 %suggest that the reader try the following exercise.

 \begin{exer}
For $Q = \R^2$ and $L = \frac{1}{2} (\dot x^2 + \dot y^2)$
the EL equations are those whose solutions are straight lines travelled at constant speed.
Rewrite $L$ in    polar coordinates $r, \theta$ and write down
the corresponding Euler-Lagrange equations,    thus deriving the equations of a straight line
in polar coordinates.   \end{exer}

 \subsection{ Reducing the least action principle}
 \label{reducedaction}
 
The curves competing in   the  least
 action principle as  we stated it are subject to    boundary conditions :  they  
 connect two    fixed points
of the configuration space $\C^3$ of   located triangles.
 Replace the  two points by two   oriented congruence classes 
 to get new boundary conditions:   the competing curves   connect
 two fixed {\it oriented congruence classes}  of configuration
 space.   If we remember  that an oriented congruence class
 is represented by  a point of shape space we 
 arrive at an  
 action principle for shape space.  
% We soon proceed to compute its Lagrangian  in shape space variables.
 
 {\it Shape space action principle.  } 
 Fix two shapes ${\bf w}_0, {\bf w}_1$  in the shape space $\R^3$. 
 Suppose that ${\bf q} (t) \in \C^3,   0, \le t \le T$  minimizes  the standard action 
 (\ref{Lag}) among 
 all curves in the space  $\C^3$ of located triangles which join the corresponding  oriented congruence classes 
 $\Sigma_0 = \pi^{-1}({\bf w}_0),  \Sigma_1 = \pi^{-1}({\bf w}_1) \subset \C^3$  in  time $T$.
 Then we will say that its projected curve $\pi({\bf q} (t)) \in \R^3$ minimizes the {\it shape space action} among
 all curves connecting the endpoints ${\bf w}_0, {\bf w}_1$ in time $T$. 

\vskip .3cm
Consider an analogous change of boundary conditions for the 
simplest action functional  in the plane, the length functional.
 Instead of minimizing the length  of curves amongst all
 curves connecting two fixed  points in the plane ,    replace these two  points
by two  disjoint circles   $\Sigma_0$  and $\Sigma_1$.
 We know that the minimizer will be a line segment
 which is perpendicular  to both $\Sigma_0$ and $\Sigma_1$ at its endpoints.
 More generally, for a   Lagrangian on $\R^n$ of   the general  form (\ref{Lag}) ,
 if we replace the fixed endpoint minimization problem with 
 the  problem of minimizing the action among all curves connecting
 two given {\it subspaces } $\Sigma_0, \Sigma_1 \subset \R^n$  then we induce a derivative
 condition at the endpoints:  namely the extremal curves, in addition to satisfying
 the Euler-Lagrange equations   must hit their targets orthogonally:     $\dot c (0) \perp \Sigma_0$ at $c(0)$
 and $\dot c(T) \perp \Sigma_1$  at $c(T)$.   We call this added  condition
 `` first variation orthogonality''.

 Returning  to our   situation where $\Sigma_0 = \pi^{-1}({\bf w}_0), \Sigma_1 = \Sigma_1 = \pi^{-1}({\bf w}_1)$
 are oriented congruence classes in $\C^3 = \R^6$. We will  interpret
   first variation orthogonality    in mechanical terms.  
   We  may sweep out all of $\Sigma_0$  by
 applying variable rigid motions  $g \in G$ to  a single point ${\bf q}_0 \in \Sigma_0$. 
 In set theory notation:  $\Sigma_0 = \{ g {\bf q}_0:  g \in G \}$.  
 Set $g = g(t)$ to form paths ${\bf q} (t) = g(t) {\bf q}_0$ in $\Sigma_0$
  and  differentiate these paths, alternately taking $g(t)$ to be a curve of  translations 
 or  a curve of rotations.  By this means  we see that the tangent space  $T_{{\bf q}_0} \Sigma_0$ to $\Sigma_0$
 at ${\bf q}_0$ is spanned by two subspaces, $\{ \dot c {\bf 1}, \dot c \in \C \}$ for 
 translations, and $\{ i \dot \theta {\bf q}_0: \dot \theta \in \R \}$: 
  \begin{eqnarray}
T_{{\bf q}_0} \Sigma_0 &  =  \text{infinitesimal rigid motions} \\
  &= \text{(translational)} + \text{(rotational)} \\
& = span_{\C} {\bf 1}  + span_{\R} (i {\bf q}) 
\end{eqnarray}
 The first  variation orthogonality condition is thus  that our extremal ${\bf q} (t)$
 be orthogonal to both the translation and rotational spaces:   $\langle  {\bf 1} , {\bf \dot q} (0) \rangle = 0$
 $\langle i \dot \theta {\bf q} (0) , {\bf \dot q} (0) \rangle = 0$.
 But as we saw in eqns (\ref{linMom}, \ref{angMom}) these orthogonality conditions are  equivalent to the assertions that 
 the linear and angular momentum are zero at ${\bf q}_0$.   (The inner product is the
 real part of the Hermitian one and  $Im(\langle {\bf q}, {\bf  \dot q} \rangle ) = Re( \langle i {\bf q}, {\bf  \dot q} \rangle$.)  
 We summarize 
 \begin{lemma}
  \label{orthogonal_lemma}
 The  curve ${\bf q}(t)$ in  $\C^3$ is orthogonal to the  oriented congruence class
 through ${\bf q}_0 = {\bf q}(0)$   if and only if its linear and angular momentum are zero at $t=0$.
 (See equations   (\ref{linMom}) and
( \ref{angMom}).)
 \end{lemma}
 Now if the curve ${\bf q}(t)$ of the lemma  is an extremal for our shape space action principle
  then it must   satisfy the EL equations which are   Newton's equations.
  Since linear and angular momentum are conserved for solutions to Newton's equations we have
  that  the linear and angular momentum are identically zero all along the curve.
  Equivalently: if an   extremal  curve is orthogonal to the $G$ orbit  $\Sigma_0$ through one of its points,
  then it is   orthogonal  to the group orbits $\Sigma_t$ through every one of its points. 
 We have established:
 \begin{proposition}  
 \label{reducedaction1}The extremals for the  shape space action principle   
 %of  minimizing the action
% among all curves connecting two given oriented congruence classes in time $T$
 are precisely those solutions to Newton's equations  whose
   linear and angular momentum are zero.  
%These are also the  solutions which are   orthogonal to the $G$-orbits. 
    \end{proposition}
 
 The proposition suggests a strategy for finding a Lagrangian $L_{shape}$  on shape space whose
 action miniimization is equivalent  the shape space action principle.  The Euler-Lagrange equations
 for $L_{shape}$ will be our desired reduced equations. 
 Break up kinetic energy into
 \begin{equation}
 K = \text{translational part} + \text{ rotational part } + \text { shape part}. \label{SaariDecomp}
 \end{equation}
 We have just agreed that the   translation and rotational part  of the kinetic energy must be   zero 
 along our shape extremals corresponding to the fact that they are orthogonal to $G$-orbits.
 Let us denote the last term, the shape term  of the kinetic energy as  $K_{shape}$.
 Thus
\begin{equation}
\label{Lshape}
L_{shape} =K_{shape}  -V
\end{equation}
 is the shape Lagrangian. It remains to express   $K_{shape}$ 
   in  terms of  shape coordinates
 $w_i$  and their time derivatives $\dot w_i$ and $V$ in terms of the $w_i$.   We find these expressions
 in the next three sections.  
 
 \subsection{Shape   kinetic energy}
 
 The  decomposition   (\ref{SaariDecomp}) applied to  veloicites  is sometimes   called the  ``Saari decomposition'':
   (\ref{Saari}, \ref{Chenciner}).  
      \begin{eqnarray}
 \dot {\bf q}  &=( \text{translational part} + \text{ rotational part } )+ \text { shape part}. \\
 & = T_q (Gq) \hskip .6cm \oplus \hskip 2.6cm  (T_q (Gq)) ^{\perp} \\
 & = \text{vertical} + \text{horizontal} 
   \label{SaariDecompVel}
 \end{eqnarray} 
%Thus (translational part) + (rotational part) $= T_q 
In the differential geometry of bundles 
   such a splitting is known as   a ```vertical- horizontal'' splitting in the theory of bundles, 
    $T_q (Gq)$ forms the ``vertical space''  and its orthogonal complement $T_q (Gq)) ^{\perp}$
   forming the ``horizontal space''. 
 This decomposition, which depends on the base point ${\bf q}$ at which the velocity is attached,   is   orthogonal  
   and leads to
 \begin{proposition} Suppose that  the
  center of mass of our located triangle is zero.  Then the Saari decomposition, eq (\ref{SaariDecomp}) above is
 $$K = \frac{1}{2}  \frac { \| {\bf P} \| ^2}{M}  
 + \frac{1}{2} \frac {J ^2}{I}  + \frac{1}{2}  \frac { \| {\bf \dot w} \| ^2}{I} $$
 where ${\bf \dot w} = \frac{d}{dt} \pi ({\bf q}(t))$and 
   $P = P({\bf \dot q}) ,   J = J({\bf q}, {\bf \dot q})$  are the  
linear and angular momenta (eq \ref{angMom}, \ref{linMom}).
In particular
\begin{equation}
\label{Kshape}
 K_{shape} =  \frac{1}{2}\frac { \| {\bf \dot w} \| ^2}{I} =  \frac{1}{2}  \frac { \| {\bf \dot w} \| ^2}{2 \sqrt{ \| {\bf w } \|}}
\end{equation}
 \end{proposition}
 
{\sc Proof.} A real  basis for the two-dimensional  translational   part  of the motion consists of 
 ${\bf 1}, i {\bf 1}$ and a real basis for the one-dimensional rotational part is  $i {\bf q}$.
The rotational part is orthogonal
 to the translational part since  the center of mass is given by     $\frac{1}{M} \langle {\bf q} , {\bf 1} \rangle$
 which we have supposed to be zero. ($M =  \langle {\bf 1} , {\bf 1} \rangle = m_1 + m_2 + m_3$ is the total mass.)
 Hence ${\bf 1}, i {\bf 1}, i {\bf q}$ is an  orthogonal basis for the vertical part,  $T_q (G q)$. 
 Normalize   to get the  orthonormal  basis  
 $$e_1, e_2, e_3 = {\bf 1}/ \sqrt{M}, i{\bf 1}/ \sqrt{M};   i {\bf q}/ \sqrt{I}$$
  for the  vertical part of the motion. 
  Let  ${\bf \dot q} \in \C^3$ be an arbitrary  vector  based at the located triangle  ${\bf q}_0 \in \C^3$ and
  expand  
  this   vector as 
  an  orthogonal direct sum to get the following  quantitative form of the Saari decomposition eq (\ref{SaariDecompVel})
  $${\bf \dot q} = \langle {\bf \dot q}, e_1 \rangle_{\R} e_1 + 
  \langle {\bf \dot q}, e_2 \rangle_{\R} e_2 +  \langle {\bf \dot q}, e_3 \rangle_{\R} e_3  + \text{(shape)}.$$
  The first three terms form 
  the vertical part of the velocity in eq (\ref{SaariDecompVel}) while the final (shape) part is, by definition, orthogonal to the first three terms
  and forms the horizontal part.
  Squaring lengths  and using the orthonormality of $e_1, e_2, e_3$ we find that
  $$\langle {\bf \dot q} , {\bf \dot q}  \rangle = |P|^2/M + J^2/ I + \text{shape}^2.$$

  It remains to show that $ |\text{shape}|^2  =  \frac{ \| \dot {\bf w} \|^2}{I} $. 
  In other words, we need to show that 
  %For this purpose, we may restrict to ${\bf \dot q}$ which are perpindicular
 % to the rotational and translational parts.  In other words, we are to
  %show that
  \begin{equation}
  \label{toshow}
   \| {\bf \dot w} \|^2 = \|{\bf q}|^2 \| {\bf \dot q}\|^2  \text{ if   } P( {\bf \dot q})= 0 ,    J ( {\bf  q}, {\bf \dot q} )= 0 \text{ and } {\bf \dot w} = D \pi_{\bf q} ({\bf \dot q}) .
   \end{equation}
To this end, write out the map $\pi^{rot}$ in real coordinates, using  $Z_j = x_j+ iy_j$,
${\bf q} = (Z_1, Z_2) = (x_1, y_1, x_2, y_2).$
We have $\pi^{rot}((x_1, y_1, x_2, y_2) = (\frac{1}{2}(x_1 ^2 + y_1 ^2 - x_2 ^2 - y_2 ^2, x_1 x_2 + y_1 y_2, x_2 y_1 - x_1 y_2).$
Compute the Jacobian: 
$$D \pi^{rot}_{\bf q}   
 = \left(
\begin{array}{cccc}
 x_1  & y_1 &-x_2 &-y_2 \\
x_2  & y_2 &x_1 & y_1 \\
-y_2  & x_2 &y_1 &-x_1  
\end{array}
\right).
$$
and set
$$L = D\pi^{rot}_{\bf q}.$$
The three rows of $L$  are orthogonal and  each has length $\|{\bf q}\|^2$.
It follows that
$$L L^T= \|{\bf q}\|^2 Id.$$
 Now the kernel of $L$ is  spanned by $e_1, e_2, e_3$  
since $\pi$ is invariant under rotations and translations.
Thus the image of $L^T$ is   the orthogonal complement to $e_1, e_2, e_3$ -- the subspace ``(shape)'' above.
Consequently any vector   ${\bf \dot q}$ of the form required in eq (\ref{toshow})
can be written   ${\bf \dot q} = L^T {\bf v}$ for some ${\bf v} \in \R^3$. 
Thus:  
\begin{eqnarray}
\| {\bf \dot q}  \|^2 & = \langle  L^T {\bf v},   L^T {\bf v} \rangle \\
& = \langle   {\bf v},  L  L^T {\bf v} \rangle \\
& = \langle    {\bf v},   \|{\bf q} \|^2 {\bf v} \rangle \\
& =   \|{\bf q} \|^2   \|{\bf v} \|^2 
\end{eqnarray}
And ${\bf \dot w} = L {\bf \dot q}$
so  that ${\bf \dot w} = L L^T {\bf v} =  \|{\bf q} \|^2 {\bf v}$.
We get that $\| {\bf \dot w}  \|^2 = \| {\bf q } \|^4  \|{\bf v} \|^2 = \|{\bf q} \|^2 \| {\bf \dot q}  \|^2 $.
Thus $\| {\bf \dot q}  \|^2 = \| {\bf \dot w}  \|^2 / \|{\bf q} \|^2$. 
Use $I = \|{\bf q} \|^2$ and $I = 2 \sqrt{\| {\bf w} \|}$ (property G of theorem \ref{main}).

\section{Shape Space metric}
\label{shapemetricsection}

\begin{definition}  The shape space metric is the twice the shape space kinetic energy $K_{shape}$, 
 viewed as a Riemannian metric on shape space, so a bilinear quadratic form on velocities depending smoothly
 on ${\bf w} \in \R^3$.  
 \end{definition} 
  We have  seen in the previous proposition that the shape space metric is given by
 \begin{equation}
 ds_{shape} ^2 =  \frac{ dw_1 ^2 + dw_2 ^2 + dw_3 ^2}{2 \sqrt{ w_1 ^2 + w_2 ^2 + w_3 ^2}}
 \label{shapemetric}
 \end{equation}
 Like any Riemannian metric, this  metric induces a distance function on shape space.
 We recall how this distance function is defined. First
 define the {\it length} $\ell$ of a path in shape space to be 
 $\ell(c) = \int_c ds_{shape} := \int_a ^b  \sqrt{2K_{shape}}dt$. 
 Now define the distance between two points as the infimum of the lengths
 of all paths that join the two points. 
 In other words,  the shape space length  is the action relative to the Lagrangian $ \sqrt{2K_{shape}}$
 and the shape space distance   between two points  is realized by an action-minimizing
 curve which joins the points.  We call such a minimizer a  {\it geodesic}.  
 
 Reparameterizing a curve does not change its length.  When we 
 parameterized a curve  by a constant
 multiple of arclength then we are insisting that  $K_{shape}$ is constant along the curve. 
 By a well-known argument involving the Cauchy-Schwartz inequality, the length minimizing curves 
 which are so parameterized are precisely the curves which minimize $\int K_{shape} dt$.
 Now the shape space action principle holds
for $K$ in place of $K-V$, and the corresponding reduced Lagrangian
is $K_{shape}$.  The geodesics for $K$ are straight lines in $\C^3$.  
Putting together these observations we have proved  the assertions of the first two sentences of:  
\begin{theorem}
\label{distancethm} The distance function defined by the Riemannian metric  agrees with the shape space distance of eq (\ref{dist}) . 
Its geodesics are the projections by $\pi: \C^3 \to  \R^3$
of horizontal lines in $\C^3$.    Each  plane $\Pi:  A w_1 + Bw_2 + Cw_3 = 0$ through the origin
is  totally geodesic: a geodesic which  starts on $\Pi$  initially
tangent to $\Pi$,  lies completely in the plane $\Pi$. The restriction of the
Riemannian metric to such a plane $\Pi$, when expressed in standard Euclidean  polar coordinates $(r, \theta)$
on that plane, has the form
$$ds^2 _{shape} |_{\Pi}  = dr^2 + \frac{1}{4} r^2 d \theta ^2.$$ 
\end{theorem} 
In order to finish the proof of this theorem, let   $\ell (t) = {\bf q} + t {\bf v}$ 
be a horizontal line passing through the point ${\bf q} \in \C^2 _0 \subset \C^3, {\bf q} \ne 0$
with horizontal tangent vector ${\bf v}$.  There are two possibilities for ${\bf v}$: a multiple of
${\bf q}$ or linearly independent of ${\bf q}$.  In the first case, we may
assume that ${\bf v} =  {\bf q}$ is the radial vector. (Note that this vector is   horizontal.) Then  
$\ell$ is a radial line and $\pi(\ell)$ is the ray  connecting  the triple collision point $0$ to ${\bf w} = \pi({\bf q}) $ (traversed twice).
The distance from $0$ to  ${\bf w}$  along this ray is the radial variable
$$r = dist(0, {\bf w}) = \|{ \bf q}\| = \sqrt{I} = \sqrt{2 \|{\bf w} \|}.$$  In the second case  
 ${\bf q}$ and   ${\bf v}$      span  a real horizontal two-plane $P$ in $\C^3$ which passes through $0$
and  contains the  line $\ell$. One computes that the  projection $\Pi:= \pi(P) \subset \R^3$ is  a plane (rel.
the coordinates $w_i$)  passing through $0$. However
  the    projection $\pi(\ell)$ is not a line (relative to the   linear coordinates  $w_i$)!

We can understand the geodesic $\pi (\ell)$ in shape space   by understanding the restriction 
$ds^2 _{shape} |_{\Pi}$ of the shape space metric
to  the plane $\pi(P)$.   Here is what we know so far about this metric.
The   radial lines  are geodesics.  The distance along such a radial geodesic 
from the triple collision point  $0$ to a random point ${\bf w} \in  \pi(P)$  is $r$ as given above.
To dilate the metric by a factor $t > 0$ we multiply ${\bf w} \in \R^3$
by $t^2$, since ${\bf q} \mapsto t {\bf q}$ corresponds to   ${\bf w} \mapsto t^2 {\bf w}$  under
$\pi$.    Finally,   the metric on $\Pi$  is rotationally symmetric, since
the   expression (\ref{shapemetric}) is rotationally invariant.
From all of this information we deduce  that the restricted  metric 
as the form:   
\begin{equation}
\label{conical}
ds ^2 _{\pi (P)} = dr^2 + c^2 r^2 d \theta ^2 
\end{equation}
where    $(r, \theta)$ are polar coordinates on the plane and $c$ is a  constant. 
It remains to show that  $c = 1/2$.   With this in mind,
consider the   circle
$r =1$ in the plane $\Pi$.  Its circumference computed from
the formula (eq \ref{conical}) is  $2 \pi c$.  
But we can also compute its length by working up 
on $P \subset \C^3$.    Take ${\bf q}$ and ${\bf v}$ to both
be  unit length and orthogonal, so an \on basis for $P$.    Then the corresponding  horizontal circle on $P$ 
is   $\cos (s) {\bf q} + \sin{s} {\bf v}$,  $0 \le s \le 2 \pi$.
But  $\pi({\bf q}) = \pi(-{\bf q})$, since $-{\bf q} = e^{i \pi}{\bf q}$,  so that the projection of this circle closes up once we have gone half
way around, from $s = 0$ to 
$s = \pi$. Thus the  projected  circle on $\pi(P) = \Pi$ has the length of half a unit circle, or   $\pi$.  
Comparing lengths we see that   $c = 1/2$.

\begin{figure}[h]
\scalebox{0.4}{\includegraphics{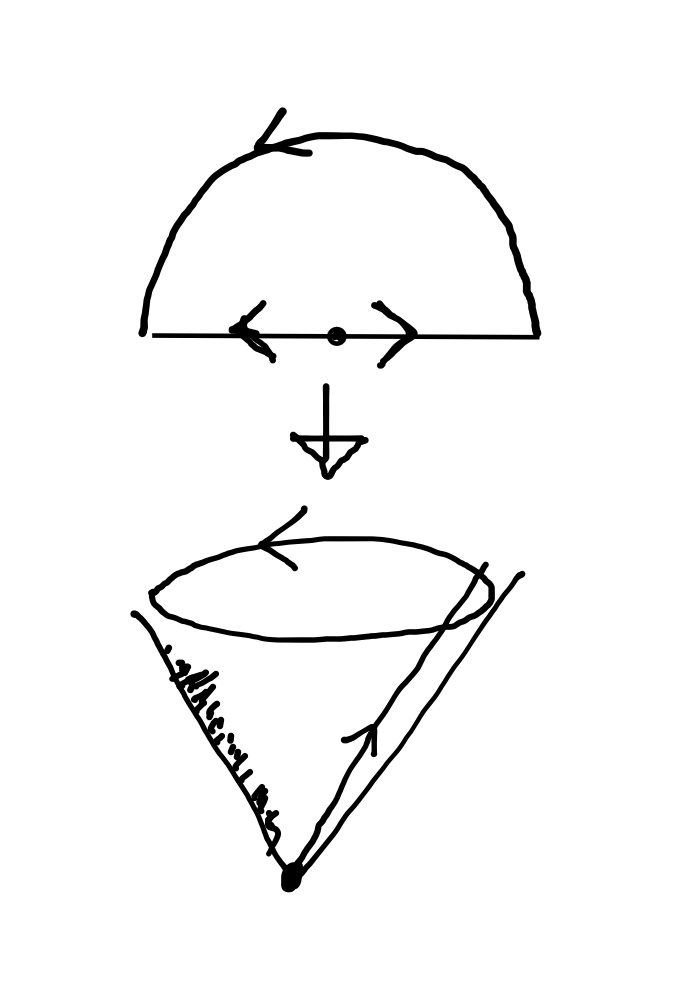}}
\caption{Folding a half-sheet of paper makes the desired cone.}   \label{conefig}
\end{figure}

Any  metric of form (\ref{conical}) is that of a cone.  We can form our $c= 1/2$  cone by taking a sheet of paper
and marking the midpoint of one edge to be the cone point. See figure \ref{conefig}. 
 Fold that edge up
so the two halves touch eachother and we have a paper model of the required cone.
Note that the circle of radius $r$ about the cone point has circumference $\pi r$,  
as required.

\section{Potential on shape space.}

In order  to express the   potential  (equation (\ref{potential})) 
in terms of the shape coordinates $w_i$ it suffices
to  express the distance $r_{ij}$ between bodies $i$ and $j$ in
terms of the $w_i$'s.    Here is the basic
geometric fact that makes this computation possible:
\begin{lemma}
Let $d_{ij}$ denote the shape-space distance from the
point ${\bf w}$ to the ij binary collision ray.
Then 
$$\sqrt{\mu_{ij}} r_{ij} =  d_{ij}$$
where $\mu_{ij}  = m_i m_j/(m_i + m_j )$. 
\end{lemma}

  Let  
${\bf b}_{12} , {\bf b}_{23}, {\bf b}_{31}$ be the  unit vectors
in shape space whose postive span defines the corresponding binary
collision ray.  For example 
${\bf b}_{12}$ represents the  binary collision $r_{12} = 0$ etc.  
The dot product and norm in  the following lemma  
 are the standard dot product and
norm on $\R^3$.  
 \begin{lemma} 
The distance $r_{ij}$ between body $i$ and $j$ is given by 
\begin{equation}
\label{sidelengths}
 r_{ij} ^2  =  {{m_i+  m_j} \over {m_i m_j}} (  \|{\bf w} \| - {\bf w} \cdot {\bf b}_{ij} ) 
\end{equation}
\end{lemma}
 
\medskip
{\bf  Proof of Lemma 2.}
We will just do the case $i,j = 1,2$.
Let
${\bf q} = (q_1, q_2, q_3) \in \C^2 _0 $ be a centered located triangle
whose projection to  shape space is 
${\bf w}$.  Let masses 1 and 2 move towards each other
along the line segment which joins them until they collide.
Make the motion linear  and such that their   center of mass remains unchanged.
Keep mass 3 fixed, so that the second Jacobi vector remains also constant.  
Described in
Jacobi coordinates this motion is given by $((1-t)Z_1, Z_2), 0 \le t \le 1$.
The curve just described is a horizontal line segment and so realizes the
shape space distance  between its endpoints in shape space $\R^3$.
One  endpoint of the segement is the initial shape ${\bf w}$. The other endpoint lies on the $12$   binary collision
ray.

When we view this line segment as a  moving triangle in the  plane, vertices  1 and 2 sweep out  
their  entire edge $[q_1, q_2]$, meeting somewhere in the middle (at their  common center of mass) while  vertex 3   remains fixed.  
From this perspective
the following is no surprise.
\begin{exer}  Compute the length of this
line segment relative to the mass metric to 
be $\sqrt{\mu_{12}} r_{12}$
where $\mu_{12} = m_1 m_2/(m_1 + m_2 )$. 
\end{exer} 

Our  Jacobi coordinate description of the 
line segment used in the proof of Lemma 2 shows that  the segment hits  the 12 collision locus $Z_1 = 0$
orthogonally at $t=1$.
Thus this segment represents a   horizontal geodesic which  minimizes
the distance from  the shape space point ${\bf w}$ to the binary collision ray associated to 
${\bf b}_{12}$. 
It follows that  the distance between ${\bf w}$ and that   collision ray is 
\begin{equation}
d_{12} = \sqrt{ \mu_{12} } r_{12}
\label{d12}
\end{equation} 
QED

\begin{figure}[h]
\scalebox{0.4}{\includegraphics{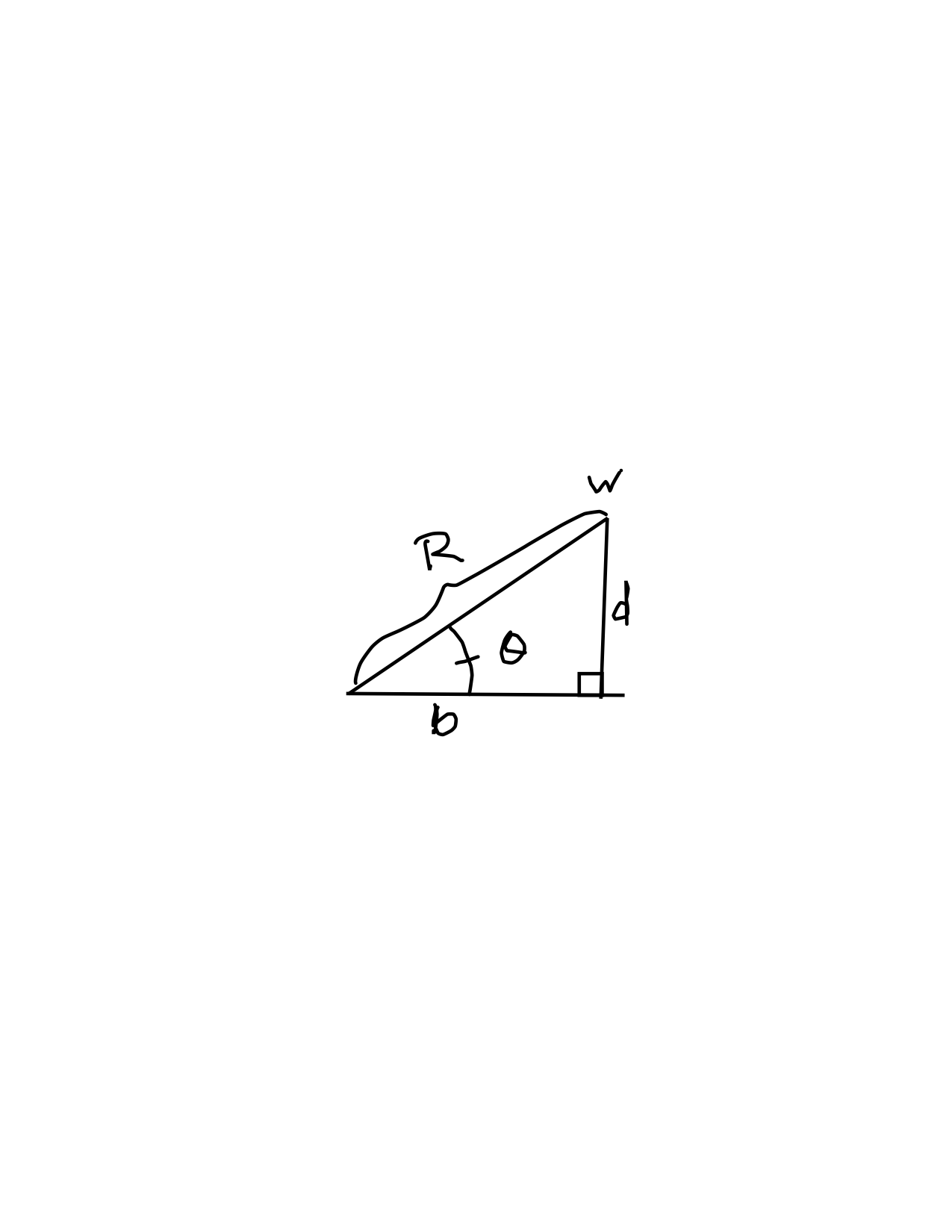}}
\caption{A Euclidean triangle allows us to compute $d =d_{12}$.}   \label{ConeTriangle}
\end{figure}

{\bf Proof of Lemma 3.} We  use the conical structure of the metric ( eq (\ref{conical}), and Theorem \ref{distancethm}) to compute
$d_{12}$ a different way.  
Two linearly independent vectors ${\bf w}, {\bf b} \in \R^3$ define two
rays which span a plane $\Pi$ through the origin. The metric on this
plane is of the form ((\ref{conical})  with $c =1/2$.  This factor of $1/2$
implies that   the shape space angle $\theta$ between the two rays is exactly  {\it half} of the
their Euclidean angle.  In other words:
$${\bf w} \cdot {\bf b} = \| {\bf w} \| cos (\psi) \; ; \psi = 2 \theta$$
The geometry of any two dimensional conical metric (\ref{conical}) is locally Euclidean.
It follows that we can compute the shape space distance $d_{12}$  between ${\bf w}$
and the ray spanned by ${\bf b}_{12}$ using   standard trigonometry
as indicated in figure \ref{ConeTriangle}.

From the figure
$d_{12} = r\sin (\theta)$
so that
$$d_{12} ^2 = r^2 \sin^2 (\theta)$$
From ${\bf w} \cdot {\bf b}_{12} = \| {\bf w}\| \cos(\psi)$,
$\psi = 2 \theta$ and $2 \| {\bf w}\|  = r^2$ we compute that
$\| {\bf w}\| - {\bf w} \cdot {\bf b}_{12} = r^2 \sin^2 \theta$
or 
\begin{equation}
d_{12}^2 =  \| {\bf w}\| - {\bf w} \cdot {\bf b}_{12} 
\label{d12squared}
\end{equation}
Formula (\ref{sidelengths}) in the case $i,j= 1,2$ now follows immediately from this last equation and
eq (\ref{d12}).

\section{Reduced equations of Motion.}

I now have written  both the shape space  kinetic energy   (eq.( \ref{Kshape}))
and
the potential energy  (eq \ref{potential}) in terms of shape space variables $w_i$.
Consequently we have the shape space   Lagrangian: 
\begin{equation}
L_{shape} =  \frac{1}{2}  \frac{ \dot w_1 ^2 + \dot w_2 ^2 + \dot w_3 ^2}{2 \sqrt{ w_1 ^2 + w_2 ^2 + w_3 ^2}} +  \frac{c_{12}}{d_{12}} +  \frac{c_{23}}{d_{23}} + \frac{c_{13}}{d_{13}} 
\label{Lshape2}
\end{equation}
Here $c_{ij} = m_i m_j \sqrt{\mu_{ij}} = (m_i m_j)^{3/2}/\sqrt{m_i + m_j}$ and $d_{ij}$ is given by formula (\ref{d12squared}) 
and is the distance betweeen the   shape space point $(w_1, w_2, w_3)$ and the $ij$ binary collision ray.  
From this expression for $L_{shape}$  we can immediately compute the equations of motion which are the Euler-Lagrange equations (\ref{EL}) for this Lagrangian:
\begin{equation}
\frac{d}{dt}( \dd{L_{shape}}{\dot w_i} )= \dd{L_{shape}}{w_i}, i =1,2,3
\label{ELshape}
\end{equation}

\section{Infinitely Many Syzygies}  The shape space  Lagrangian (eq \ref{Lshape2}) is that of a point mass moving in $\R^3$
(endowed with   metric (\ref{shapemetric}) )
subject to  the attractive force generated by the pull of the three
binary collision rays.   These rays   lie in the
collinear plane $w_3 =0$.  Consequently the point ${\bf w}$ is always attracted toward
the collinear plane.   This  suggests that the shape point ${\bf w}$ must oscillate
 back and forth crossing that plane infinitely often.  Indeed :  
\begin{theorem} [See \cite{InfinitelyMany}] 
\label{syzygies} If  a solution with negative energy and zero angular momentum   does not
begin and end in triple collision then it must cross  the collinear plane $w_3 = 0$ infintely often.
\end{theorem}
{\sc  Sketch of proof of theorem }\ref{syzygies}. 
 The  heuristic physical description of the reduced dynamics  described just before stating
 the theorem  led 
 us to discover 
a differential equation of the form $\frac{d}{dt}  (f \frac{d}{dt}z ) = - g z$
for a normalized height variable $z= w_3/\tilde I$
where $\tilde I$ is the  homogeneous quadratic positive  function
of the $w_i$. ($\tilde I$ is the function  $I$ with all the masses artificially assumed   equal.) 
Here  $f$  is a positive  function on shape space and  $g$ is a non-negative  function of the $w_i$
and $\dot w_i$ which is   positive
away from the   Lagrange homothety solution.
The result  follows  from this differential equation by a kind of Sturm-Liousville argument.

   \section{Finale}   We end this article with  another   
theorem whose conception and proof was made possible
through the  shape space formulation of the three-body problem 
 
 \begin{theorem}[See \cite{Remarkable}]
\label{eight}
There is a periodic solution to the equal mass zero angular momentum three-body problem
in which all three masses chase each other around the same figure-eight shaped curve.
\end{theorem}

%\end{document}
 
\begin{figure}[h]
\scalebox{0.5}{\includegraphics{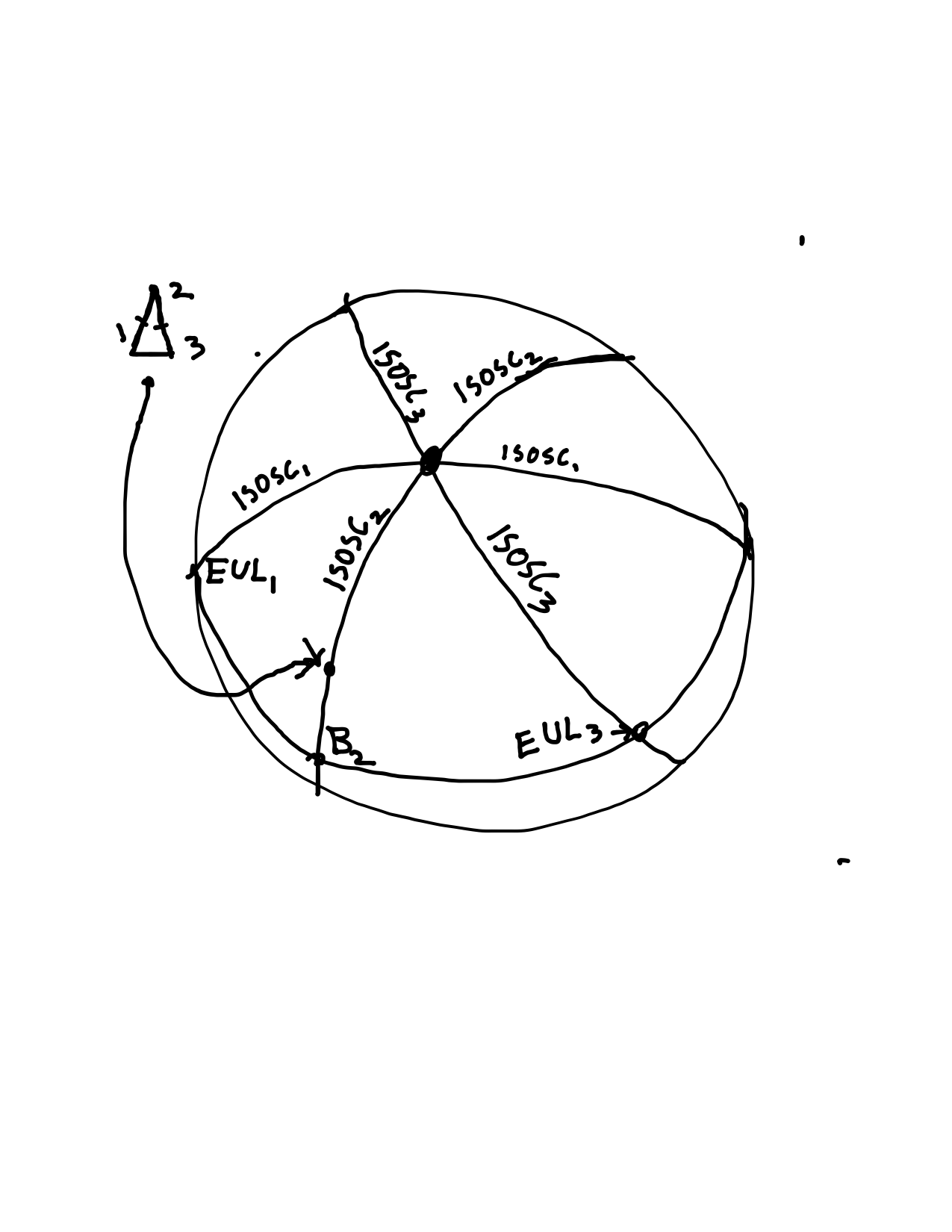}}
\caption{The shape sphere with the  Isosceles great circles marked.}  
\label{spherewIsosc}
\end{figure}

{\bf Sketch of proof of theorem \ref{eight}.}   The space   of isosceles triangles forms three great circles 
on the shape sphere,  each circle passing through the   two Lagrange points which are represented as the North and South Poles.
See   figure \ref{spherewIsosc}.  

Label these circles $ISOSC_1, ISOSC_2, ISOSC_3$  according to the   mass  which forms the isosceles triangle's vertex.
For example   $ISOSC_1$  is the circle consisting of those triangles for which $r_{12} = r_{13}$.
Each circle $ISOSC_i$  intersects the collinear plane $w_3 = 0$  in two points,  an Euler point
$EUL_i$ which represents  the
configuration of Euler's solution, and  a binary collision point $B_i$ (our old ${\bf b}_{jk}$). 
For example $EUL_1 \in ISOSC_1$ represents those degenerate triangles for which mass   1 lies at the midpoint of the segment formed
by masses 2 and 3.
%, and the collision point $B_1$ represents those triangles for which   $r_{23} = 0$.  
Considered in  shape space, the circle $ISOSC_i$ becomes a plane $\Pi_i$ and    
  $EUL_i$ and $B_i$ becomes a ray lying in that plane.   

We consider the problem of  minimizing the shape space action $\int_c L_{shape} dt$ among all paths connecting $EUL_1$ to $ISOSC_2$ 
in a fixed time $T$.  The difficult part of the proof is showing that this minimimum  $\gamma_*$ actually  exists and is collision-free.
Once established, we know by the first variation orthogonality that $\gamma_*$ must hit the Euler ray   $EUL_1$ orthogonally, and the 
isosceles plane $ISOSC_2$ orthogonality.
This orthogonality  allows us to  continue the solution arc $\gamma_*$ by reflection.  The equal mass condition
 insures that these reflections take solutions to solutions.   

For example, the fact that $m_1 =m_3$ implies that  permutation of masses 1 and 3 is a  linear isometric
involution  $\sigma_2: (q_1,q_2,q_3) \mapsto (q_3, q_2, q_1)$
of $\C^3$ (isometric relative to the mass metric) which  preserves the potential and hence maps solutions to Newton's equations to 
solutions to Newton's equations.   The effect of $\sigma_2$ on the shape sphere
is a half-twist about the line connecting $EUL_2$ and $B_2$.   If we compose $\sigma_2$ with
the operation of  reflection about the symmetry axis of 
isosceles triangle  $\gamma_* (T)$ we obtain
an action-preserving 
isometric involution  $R_2$ whose the effect 
on shape space is that of   reflection about  the plane $ISOSC_2$.
Thus $R_2 (\gamma_*)$ is a solution arc  which connects
$ISOSC_2$ to $EUL_3$ and whose derivative matches up with $\gamma_*$. It therefore agrees
with the continuation of the solution $\gamma_*$, up to a time translation. We now have
a solution on the interval $[0,2T]$ joining $EUL_1$ to $EUL_3$. 
Apply the half-twist  $\sigma_3$ about the line connecting
$EUL_3$ and $B_3$  to this  continued solution  
to obtain a solution arc on the interval $[0,4T]$ connecting $EUL_1$ to $EUL_2$,
and crossing $ISOSC_2, EUL_3, ISOSC_1$ in that order along the way. 
Continue to march around the  sphere, applying half-twists or reflections as
needed   I  obtain  twelve solution  arcs  each congruent to the original arc $\gamma_*$ and   whose derivatives
all  match up. The result is a   smooth   solution to the three-body equation  which is {\it periodic in shape space} with period $12T$.   
 With some extra work, we can show that the solution is actually periodic of period $12T$ in the
 original inertial space -- i.e the plane, and that this   curve $\gamma$  is a 
 {\it choreograpy}.   The curve $\gamma$ is the figure eight solution. 
 
 \begin{definition}  An N-body choreography of period $T$ is a solution to N-body problem
 which has the particula form
 ${\bf q} (t) = (q_1 (t), q_2 (t), \ldots , q_N (t))$  where $q_i (t) = s (t - i T/N) $, $i = 2, 3, \ldots , N$
 for some fixed  $T$ periodic curve  $s(t)$   in the plane (or space). 
 \end{definition}

\begin{figure}[h]
\scalebox{0.5}{\includegraphics{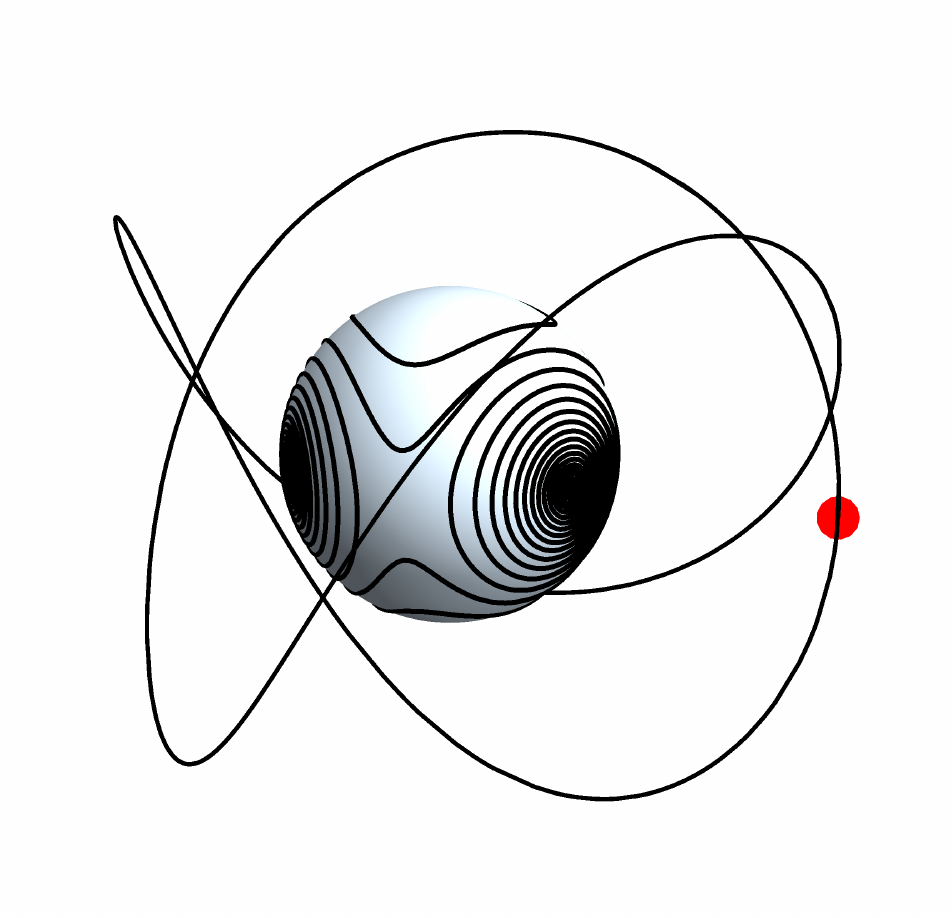}}
\caption{The figure eight as it appears in shape space.}  
\label{eightinshape}
\end{figure}

This curve  is the figure eight, viewed in shape space.  See figure \ref{eightinshape}.

{\bf History.} The figure eight solution curve was discovered   in 1994 using the principle of least action by C. Moore \cite{Moore}. 
Chenciner and myself rediscovered the solution using the shape space   least action principle in 2000. Our methods yielded 
  a rigorous existence proof.
   
 \section{Summary}

The perspective  of  shape space  reduces many features of   the  Newtonian
three-body problem  to an investigation of curves on the shape sphere  or even the shape space. 
These spaces have dimension 2 and 3 compared to the 12 dimensions of the original planar three-body phase space:
 6   for the  vertices of the moving triangle and 6     for their velocities.   Our minds are not  nearly  as well attuned to  
12 dimensional visualization  as they are to  3 dimensional visualization.   The reduction from 12 to 3  
allows us  hope of accurately and successfully  applying the  creative  powers of our geometric
visualization to better understand the three-body problem.    
% We hope we have convinced the reader of the utility of this reduction in dimension.

 \section{Acknowledgements.} I would like to thank Alain Chenciner, Alain Albouy, Carles Simo
 and Rick Moeckel for many conversations over the years. 
  I would particularly like to
% thank Alain Albouy for teaching me how to view Heron's formula as defining a light cone
  Piet Hut for the use of  figure \ref{Hut} and Rick Moeckel for the use of 
 figure \ref{eightinshape}.

\end{document}